\documentclass[12pt]{amsart}
\usepackage[mathscr]{eucal}
\usepackage{amssymb,amsmath}
\usepackage{amsfonts}
\usepackage{a4}

\newtheorem{theorem}{Theorem}
\newtheorem{proposition}{Proposition}
\newtheorem{lemma}{Lemma}
\newtheorem{corollary}{Corollary}

\begin{document}
\title{Noncommutative geometry through monoidal categories}

\author[Tomasz Maszczyk]{Tomasz Maszczyk\dag }
\address{Institute of Mathematics\\
Polish Academy of Sciences\\
Sniadeckich 8\newline 00--956 Warszawa, Poland\\
\newline Institute of Mathematics\\
University of Warsaw\\ Banacha 2\newline 02--097 Warszawa, Poland}
\email{maszczyk@mimuw.edu.pl}

\thanks{\dag The author was partially supported by KBN grants 1P03A 036 26 and 115/E-
343/SPB/6.PR UE/DIE 50/2005-2008.}
\subjclass{14A22, 16S38, 16W30,
16E40}

\begin{abstract} After introducing a noncommutative counterpart of
commutative algebraic geometry based on monoidal categories of
quasi-coherent sheaves we show that various constructions in
noncommutative geometry (e.g. Morita equivalences, Hopf-Galois
extensions) can be given geometric meaning extending their
geometric interpretations in the commutative case. On the other
hand, we show that some constructions in commutative geometry
(e.g. faithfully flat descent theory, principal fibrations,
equivariant and infinitesimal geometry) can be interpreted as
noncommutative geometric constructions applied to commutative
objects. For such generalized geometry we define global invariants
constructing cyclic objects from which we derive Hochschild,
cyclic and periodic cyclic homology (with coefficients) in the
standard way.

\end{abstract}

\maketitle

\tableofcontents

\section{Introduction} Abelian categories  as a replacement for
spaces (schemes) can be justified by the following reconstruction
theorem.

\paragraph{\textbf{Theorem.}} (P. Gabriel for noetherian schemes (\cite{gab},
Ch. VI); A.L. Rosenberg in quasicompact case (\cite{ros}); and in
general case ([115])) {\em Every scheme X can be reconstructed
from the abelian category QcohX with the distinguished object
${\mathcal O}_{X}$ uniquely up to an isomorphism of schemes.}

Morphisms between schemes are encoded on the level of
quasi-coherent sheaves as pairs of adjoint functors (the direct
image and the inverse image as its left adjoint), in a way
resembling geometric morphisms among topoi \cite{macmoe}.

The idea of a noncommutative algebraic geometry, based on abelian
categories or their generalizations (triangulated categories,
dg-categories and $A_{\infty}$-categories) \cite{mur,
gol,art,ver,ros,konros,orl} is derived from the following
observation. The category of modules makes sense for any
associative, not necessarily commutative, ring. Therefore
arbitrary (with some working restrongions) abelian (or
triangulated, dg, $A_{\infty}$) categories should be regarded as
categories of quasi-coherent sheaves (or complexes of sheaves) on,
possibly non-affine, non-commutative ``schemes". This theory
develops in close relation with representation theory
\cite{hod,lunros'''}.

However, in this approach one important point from commutative
geometry is lost. Classical algebraic geometry is based on
polynomials. They describe varieties and morphisms between them.
Composition of morphisms is defined by substitution of polynomials
into polynomials. The natural environment for polynomials are
symmetric monoidal categories, and categories of quasi-coherent
sheaves are such. Polynomial  substitutions produce (co)monoidal
functors between these monoidal categories. Lack of monoidal
structures is the main drawback of module categories over
noncommutative rings. Although one can derive from a module
category its monoidal category of bimodules regarded as
endofunctors \cite{wat,eil}, in general there is no way to
transport them along module-theoretic geometric morphisms, and if
it is accidentally possible, the result is different from the
result obtained for symmetric bimodules (over a commutative ring)
regarded simply as modules.

One could argue that modules are important because of
representation theory. But group algebras and enveloping algebras
of Lie algebras are augmented algebras and modules over them can
be regarded as bimodules (symmetric over a ground field) with the
second side defined by means of the augmentation. Note that as
such they can be used as coefficients of Hochschild (co)homology
computing group and Lie algebra (co)homology.

Many natural constructions on noncommutative rings (or algebras)
produce bimodules (algebras, ideals, universal differentials).
Explicit natural modules for such rings, different from natural
bimodules with one side forgotten, in general are not known.

The aim of the present paper is to persuade monoidal categories as
models of quasicoherent sheaves on noncommutative schemes. This
approach is justified by the monoidal version of the
reconstruction theorem due to Balmer \cite{bal}. From this
perspective, algebras and coalgebras are not primary objects but
artifacts of geometric morphisms between noncommutative schemes.
Instead of thinking of classical spaces as of commutative
algebras, we think of abelian symmetric monoidal categories. Since
even commutative algebras admit non-symmetric bimodules, this
provides some room to consider non-classical (non-local) effects
even for classical spaces. We show that in this framework one can
study global and infinitesimal structures of a noncommutative
scheme. We compare purely geometric constructions (i.e. these
which use only some geometric morphisms on the purely categorial
level) and purely algebraic constructions (i.e. these which use
homomorphisms of some algebraic structures).

In the global picture we prove theorem about equivalence of flat
covers in the category of noncommutative affine schemes and
noncommutative Galois extensions. It means that descent data or
coactions, which are encoded in  comodule structures, can be
understood as geometric gluing or geometric quotiening  by
symmetries.

In the infinitesimal picture we establish a noncommutative duality
between \emph{infinitesimals} and \emph{differential operators}
(well known in the classical situation) realized by passing to the
opposite category. To achieve this we prove that infinitesimals
and differential operators arise as specializations of two dual
categorial constructions.

Finally, we construct global invariants of our noncommutative
schemes. They are Hochschild, cyclic and periodic cyclic homology
derived from cyclic objects. Our construction allows to introduce
coefficients into the theory, which are noncommutative analogs of
sheaves with integrable connection from the theory of the DeRham
cohomology. In a sense, the respective ``integrability condition"
in terms of some braiding is as general as possible, because it is
derived from the very structure of the cyclic object, in opposite
to other approaches where it is based on some ideas from category
theory (comonads and distributivity laws in \cite{bost} or
symmetric monoidal categories and cocartesian objects in
\cite{kal}) producing some cyclic objects. We compare different
types of diagrams standing behind our construction and
constructions based on these categorial ideas.

The present paper is a part of some kind ``noncommutative EGA in a
nutshell", \textit{tout proportion garde}, whose further topics
will appear in subsequent papers.

\section{Noncommutative schemes as monoidal categories}

\subsection{Category of noncommutative schemes}
\paragraph{\textbf{Definition.}} We define the
category $\mathfrak{Sch}$ of (noncommutative) schemes as follows.
Objects of $\mathfrak{Sch}$, usually denoted by $X$, are abelian
monoidal categories, usually denoted by $({\rm {\rm Qcoh}}(X),
\otimes_{X}, \mathcal{O}_{X})$. Morphisms $f: X\rightarrow Y$ are
isoclasses of pairs $(f_{*}, \mathcal{O}_{Y}\rightarrow
f_{*}\mathcal{O}_{X})$, where $f_{*}$ is an additive monoidal
functor $f_{*}: {\rm Qcoh}(X)\rightarrow {\rm Qcoh}(Y)$ having the
left adjoint $f^{*}$, and $\mathcal{O}_{Y}\rightarrow
f_{*}\mathcal{O}_{X}$ is a morphism in the category ${\rm
Qcoh}(Y)$, with natural composition.

\vspace{3mm}\paragraph{\textbf{Remark.}} Morphisms $f:
X\rightarrow Y$ in $\mathfrak{Sch}$ can be equivalently defined as
isoclasses of pairs $(f^{*}, f^{*}\mathcal{O}_{Y}\rightarrow
\mathcal{O}_{X})$, where $f^{*}$ is an additive comonoidal functor
$f^{*}: {\rm Qcoh}(Y)\rightarrow {\rm Qcoh}(X)$ having the right
adjoint $f_{*}$, and $f^{*}\mathcal{O}_{Y}\rightarrow
\mathcal{O}_{X}$ is a morphism in the category ${\rm Qcoh}(X)$,
with natural composition.

\vspace{3mm}\paragraph{\textbf{Example 1. (Commutative schemes)}}
With every commutative scheme $X$ one can associate its abelian
category ${\rm Qcoh}(X)$ of complexes of quasicoherent sheaves,
with the distinguished structural sheaf $\mathcal{O}_{X}$. The
tensor product $\otimes_{X}$ of $\mathcal{O}_{X}$-modules makes
${\rm Qcoh}(X)$ a monoidal category with $\mathcal{O}_{X}$ as the
unit object. With every morphism of commutative schemes $f:
X\rightarrow Y$ one can associate the additive monoidal (direct
image) functor $f_{*}: {\rm Qcoh}(X)\rightarrow {\rm Qcoh}(Y)$,
which has the left adjoint $f^{*}$ (the inverse image functor),
and a morphism $\mathcal{O}_{Y}\rightarrow f_{*}\mathcal{O}_{X}$
in the category ${\rm Qcoh}(Y)$.

This example has some special features. The monoidal categories
${\rm Qcoh}(X)$ are symmetric, the direct images
$f_{*}\mathcal{O}_{X}$ are commutative algebras in symmetric
monoidal categories ${\rm Qcoh}(Y)$, and finally, the inverse
image functor $f^{*}$ is strongly comonoidal.

Instead of the category of quasicoherent sheaves one can consider
the derived category of perfect complexes, with its canonical
monoidal structure. The benefit from this upgrading is the
reconstruction theorem of Balmer \cite{bal}, which provides a
construction on symmetric tensor triangulated categories with
values in locally ringed spaces, functorial with respect to all
tensor triangulated functors, reconstructing a topologically
noetherian scheme from its derived category of perfect complexes.

\vspace{3mm}\paragraph{\textbf{Example 2. (Finite flat
correspondences of commutative schemes)}} One can consider
category, whose objects are commutative schemes but morphisms $f$
from a scheme $X$ to a scheme $Y$ are defined as isoclasses of
diagrams  of the form
$$\begin{array}{ccc}
          \widetilde{X} & \stackrel{\widetilde{f}}{\longrightarrow} & Y  \\
\!\!\!\!\!\pi\downarrow &                                           &    \\
          X             &                                           &
\end{array}$$
in the category of schemes, with $\pi$ finite flat and
$\widetilde{f}$ separable and quasi-compact (this is a technical
assumption on $\widetilde{f}$ for the flat base change
isomorphism). The composition of morphisms is defined by means of
the following diagrams

$$\begin{array}{ccccc}
\widetilde{\widetilde{X}}           & \stackrel{\widetilde{\widetilde{f}}}{\longrightarrow} &  \widetilde{Y}            & \stackrel{\widetilde{g}}{\longrightarrow}   & Z\\
\!\!\!\!\!\widetilde{\pi}\downarrow &         \square                                       &  \ \ \downarrow\rho &                                             & \\
\widetilde{X}                       & \stackrel{\widetilde{f}}{\longrightarrow}             & Y                         &                                             & \\
\!\!\!\!\!\pi\downarrow             &                                                       &                           &                                             & \\
X                                   &                                                       &                           &                                             & \\
\end{array}$$
where $\square$ denotes a cartesian square.

Every such a morphism $f=(\pi, \widetilde{f})$ defines a functor
\begin{align}f_{*}:= \widetilde{f}_{*}\pi^{*}: {\rm
Qcoh}(X)\rightarrow{\rm Qcoh}(Y).
\end{align}
Since $\widetilde{f}_{*}$ is monoidal and $\pi^{*}$ is strongly
comonoidal, $f_{*}$ is monoidal as well. There is also a canonical
homomorphism
\begin{align}
{\mathcal O}_{Y}\rightarrow f_{*}{\mathcal
O}_{X}=\widetilde{f}_{*}\pi^{*}{\mathcal
O}_{X}=\widetilde{f}_{*}{\mathcal O}_{\widetilde{X}}.
\end{align}
of algebras in ${\rm Qcoh}(Y)$. To prove that $f_{*}$ has the left
adjoint it is enough to prove this property for $\pi^{*}$.

\begin{lemma}
For every $\pi: \widetilde{X}\rightarrow X$ finite flat
$\pi^{*}={\mathcal H}om_{X}((\pi_{*}{\mathcal
O}_{\widetilde{X}})^{\vee}, {\mathcal F})^{\sim}$ and has the left
adjoint $(\pi_{*}{\mathcal
O}_{\widetilde{X}})^{\vee}\otimes_{\pi_{*}{\mathcal
O}_{\widetilde{X}}}\pi_{*}(-)$.
\end{lemma}

\emph{Proof.} Notice first that being affine $\pi$ satisfies
\begin{align}
\pi^{*}{\mathcal F}=(\pi_{*}{\mathcal
O}_{\widetilde{X}}\otimes_{X}{\mathcal F})^{\sim}.
\end{align}
Since $\pi_{*}{\mathcal O}_{\widetilde{X}}$ is locally free
coherent on $X$ the latter can be rewritten as follows
\begin{align}
(\pi_{*}{\mathcal O}_{\widetilde{X}}\otimes_{X}{\mathcal
F})^{\sim}={\mathcal H}om_{X}((\pi_{*}{\mathcal
O}_{\widetilde{X}})^{\vee}, {\mathcal F})^{\sim}.
\end{align}
Here  $(\pi_{*}{\mathcal O}_{\widetilde{X}})^{\vee}$ is equipped
with the canonical contragredient  $\pi_{*}{\mathcal
O}_{\widetilde{X}}$-module structure. Using the fact that
$\pi_{*}$ is an equivalence between ${\rm Qcoh}(\widetilde{X})$
and the subcategory of ${\rm Qcoh}(X)$ of $\pi_{*}{\mathcal
O}_{\widetilde{X}}$-modules with $\pi_{*}{\mathcal
O}_{\widetilde{X}}$-linear morphisms  with the inverse
$(-)^{\sim}$ and the tensor-hom adjunction we obtain the left
adjoint of $\pi^{*}$ of the form  $(\pi_{*}{\mathcal
O}_{\widetilde{X}})^{\vee}\otimes_{\pi_{*}{\mathcal
O}_{\widetilde{X}}}\pi_{*}(-)$.   $\blacksquare$

\vspace{3mm}\paragraph{\textbf{Definition.}} Using the fact that
the ${\mathcal O}_{X}$-module $\pi_{*}{\mathcal
O}_{\widetilde{X}}$ is locally free coherent on $X$ we can dualize
the unit ${\mathcal O}_{X}\rightarrow\pi_{*}{\mathcal
O}_{\widetilde{X}}$  and the multiplication $\pi_{*}{\mathcal
O}_{\widetilde{X}}\otimes_{X}\pi_{*}{\mathcal
O}_{\widetilde{X}}\rightarrow \pi_{*}{\mathcal O}_{\widetilde{X}}
$ to define a counital coalgebra $D:=(\pi_{*}{\mathcal
O}_{\widetilde{X}})^{\vee}$. Then the algebra $\pi_{*}{\mathcal
O}_{\widetilde{X}}$ itself can be regarded as a convolution
algebra ${\mathcal H}om_{X}(D, {\mathcal O}_{X} )$ of the
coalgebra $D$, hence $\widetilde{X}={\rm Spec}_{X}({\mathcal
H}om_{X}(D, {\mathcal O}_{X} ))$.

\begin{corollary} If $f=(\pi, \widetilde{f})$ is a
finite flat correspondence as above, then
\begin{align}f_{*}=\widetilde{f}_{*}({\mathcal H}om_{X}(D, -)^{\sim})\end{align} and has the left
adjoint
\begin{align}f^{*}:=D\otimes_{\pi_{*}{\mathcal
O}_{\widetilde{X}}}\pi_{*}\widetilde{f}^{*}(-).\end{align}
Moreover, there is a canonical homomorphism of quasicoherent
algebras on $Y$
\begin{align}{\mathcal O}_{Y}\rightarrow f_{*}{\mathcal O}_{X}=\widetilde{f}_{*}{\mathcal O}_{\widetilde{X}}.\end{align}

\end{corollary}

Finally, we have the following lemma.

\begin{lemma}
If $f=(\pi, \widetilde{f})$ and $g=(\rho, \widetilde{g})$ are
finite flat correspondences as above, then
\begin{align}g_{*}f_{*}=(gf)_{*}.
\end{align}
\end{lemma}

\emph{Proof.} By the flat base change formula applied to the
cartesian square in the above composition of correspondences we
have
\begin{align}
\rho^{*}\widetilde{f}_{*}=\widetilde{\widetilde{f}}_{*}\widetilde{\pi}^{*}
\end{align}
which implies that
\begin{align}g_{*}f_{*} & =\widetilde{g}_{*}\rho^{*}\widetilde{f}_{*}\pi^{*}=\widetilde{g}_{*}\widetilde{\widetilde{f}}_{*}\widetilde{\pi}^{*}\pi^{*}=(\widetilde{g}\widetilde{\widetilde{f}})_{*}(\pi\widetilde{\pi})^{*}=
(gf)_{*}.\ \blacksquare
\end{align}

Now we are to prove that this composition of functors $g_{*}$ and
$f_{*}$ can be described by means of coalgebras which define these
functors as in Corollary 1. We need for that an easy base change
formula with affine morphisms.

\begin{lemma} For an arbitrary cartesian square of schemes
$$\begin{array}{ccccc}
\widetilde{\widetilde{X}}           & \stackrel{\widetilde{\widetilde{f}}}{\longrightarrow} &  \widetilde{Y}            \\
\!\!\!\!\!\widetilde{\pi}\downarrow &         \square                                       &  \ \ \downarrow\rho \\
\widetilde{X}                       &
\stackrel{\widetilde{f}}{\longrightarrow}             & Y
\end{array}$$
with $\rho$ affine, the natural base change transformation
\begin{align}
\widetilde{f}^{*}\rho_{*}\rightarrow
\widetilde{\pi}_{*}\widetilde{\widetilde{f^{*}}}
\end{align}
is an isomorphism.
\end{lemma}

\emph{Proof.} Since $\widetilde{f}^{*}$ is strong comonoidal it
transforms the ${\mathcal O}_{Y}$-algebra $\rho_{*}{\mathcal
O}_{\widetilde{Y}}$ into an ${\mathcal O}_{\widetilde{X}}$-algebra
$\widetilde{f}^{*}\rho_{*}{\mathcal
O}_{\widetilde{Y}}=\widetilde{\pi}_{*}{\mathcal
O}_{\widetilde{\widetilde{X}}}$ and every $\rho_{*}{\mathcal
O}_{\widetilde{Y}}$-module ${\mathcal G}$ quasicoherent on $Y$
into a $\widetilde{\pi}_{*}{\mathcal
O}_{\widetilde{\widetilde{X}}}$-module $\widetilde{f}^{*}{\mathcal
G}$  quasicoherent on $\widetilde{X}$. Let us denote by ${\mathcal
G}^{\sim}$ and $(\widetilde{f}^{*}{\mathcal G})^{\sim}$ the
corresponding quasicoherent sheaves on $\widetilde{Y}={\rm
Spec}_{Y}(\rho_{*}{\mathcal O}_{\widetilde{Y}})$ and
$\widetilde{\widetilde{X}}={\rm
Spec}_{\widetilde{X}}(\widetilde{\pi}_{*}{\mathcal
O}_{\widetilde{\widetilde{X}}})$, respectively. Then
\begin{align}
\widetilde{\widetilde{f^{*}}}({\mathcal
G}^{\sim})=(\widetilde{f}^{*}{\mathcal G})^{\sim}.
\end{align}
Applying now $\widetilde{\pi}_{*}$ and the fact that ${\mathcal
G}=\rho_{*}({\mathcal G}^{\sim})$ we obtain
\begin{align}
\widetilde{\pi}_{*}\widetilde{\widetilde{f^{*}}}({\mathcal
G}^{\sim})=\widetilde{f}^{*}\rho_{*}({\mathcal G}^{\sim}).
\end{align}
This ends the proof of the lemma, because all quasicoherent
sheaves on $\widetilde{Y}={\rm Spec}_{Y}(\rho_{*}{\mathcal
O}_{\widetilde{Y}})$ are of the form ${\mathcal G}^{\sim}$.
$\blacksquare$

Now we can prove that the composition of inverse image functors
corresponding to finite flat correspondences can be defined in
terms of coalgebras corresponding to them according to  Corollary
1, without any reference to cartesian squares of schemes.
\begin{lemma}
If $D$ and $E$ are coalgebras defining functors $f^{*}$ and
$g^{*}$, respectively, accordind to Corollary 1, the composition
$f^{*}g^{*}$ is defined by a counital coalgebra
$D\otimes_{\pi_{*}{\mathcal
O}_{\widetilde{X}}}\pi_{*}\widetilde{f}^{*}E$.
\end{lemma}

\emph{Proof.} First we have to prove that coalgebra structures on
$D$ and $E$ define a coalgebra structure on
$D\otimes_{\pi_{*}{\mathcal
O}_{\widetilde{X}}}\pi_{*}\widetilde{f}^{*}E$. It is so because
$\pi_{*}$ being affine is a strong monoidal equivalence between
${\rm Qcoh}(\widetilde{X})$ and the subcategory of ${\rm Qcoh}(X)$
consisting of $\pi_{*}{\mathcal O}_{\widetilde{X}}$-modules with
$\pi_{*}{\mathcal O}_{\widetilde{X}}$-linear morphisms and the
tensor product over $\pi_{*}{\mathcal O}_{\widetilde{X}}$. The
counit is defined as the natural composite
\begin{align}
D\otimes_{\pi_{*}{\mathcal
O}_{\widetilde{X}}}\pi_{*}\widetilde{f}^{*}E\rightarrow
D\otimes_{\pi_{*}{\mathcal
O}_{\widetilde{X}}}\pi_{*}\widetilde{f}^{*}{\mathcal
O}_{Y}=D\otimes_{\pi_{*}{\mathcal
O}_{\widetilde{X}}}\pi_{*}{\mathcal
O}_{\widetilde{X}}=D\rightarrow{\mathcal O}_{X}.
\end{align}

Now the composition. Below we apply the following facts:
\begin{align}
(-)\otimes_{\rho_{*}{\mathcal O}_{\widetilde{Y}}}(-)={\rm
coker}((-)\otimes_{Y}\rho_{*}{\mathcal
O}_{\widetilde{Y}}\otimes_{Y}(-)\rightarrow (-)\otimes_{Y}(-)),
\end{align}
$f^{*}$ is strong comonoidal and right exact, $\pi_{*}$  is
strong monoidal as above and (right) exact, $\rho$ is affine and
we use Lemma 3.
\begin{align}
f^{*}g^{*} & = D\otimes_{\pi_{*}{\mathcal
O}_{\widetilde{X}}}\pi_{*}\widetilde{f}^{*}(E\otimes_{\rho_{*}{\mathcal
O}_{\widetilde{Y}}}\rho_{*}\widetilde{g}^{*}(-))\\
           & = D\otimes_{\pi_{*}{\mathcal
O}_{\widetilde{X}}}\pi_{*}\widetilde{f}^{*}E\otimes_{\pi_{*}\widetilde{f}^{*}\rho_{*}{\mathcal
O}_{\widetilde{Y}}}\pi_{*}\widetilde{f}^{*}\rho_{*}\widetilde{g}^{*}(-))\\
           & = D\otimes_{\pi_{*}{\mathcal
O}_{\widetilde{X}}}\pi_{*}\widetilde{f}^{*}E\otimes_{\pi_{*}\widetilde{\pi}_{*}\widetilde{\widetilde{f^{*}}}{\mathcal
O}_{\widetilde{Y}}}\pi_{*}\widetilde{\pi}_{*}\widetilde{\widetilde{f^{*}}}\widetilde{g}^{*}(-))\\
           & = (D\otimes_{\pi_{*}{\mathcal
O}_{\widetilde{X}}}\pi_{*}\widetilde{f}^{*}E)\otimes_{(\pi\widetilde{\pi})_{*}{\mathcal
O}_{\widetilde{\widetilde{X}}}}(\pi\widetilde{\pi})_{*}(\widetilde{g}\widetilde{\widetilde{f}})^{*}(-)).\
\blacksquare
\end{align}
In the next example we will recognize Corollary 1 and Lemma 4 in
their noncommutative versions.

 \vspace{3mm}\paragraph{\textbf{Example 3.
(Noncommutative affine schemes)}}

\vspace{3mm}\paragraph{\textbf{Definition.}} We call
\emph{(convolution) representation} of a ring $B$ over another
ring $A$ an arbitrary ring homomorphism $B\rightarrow {\rm
Hom}_{A}(D,A)_{A}$ to the convolution ring of a given $D\in{\rm
Coalg}(A)$. Two representations $B\rightarrow {\rm
Hom}_{A}(D,A)_{A}$ and $B\rightarrow {\rm Hom}_{A}(D',A)_{A}$ are
\emph{isomorphic} if there is an isomorphism of counital
coalgebras $D\rightarrow D'$ making the diagram
$$\begin{array}{ccc}
B & \rightarrow & {\rm Hom}_{A}(D,A)_{A}\\
\parallel &     & \uparrow\\
B & \rightarrow & {\rm Hom}_{A}(D',A)_{A}
\end{array}$$
commutative.

 \vspace{3mm}\paragraph{\textbf{Example.}} Let  $P$ be a
finitely generated projective right $A$-module and $D=P^{*}\otimes
P$ be a coalgebra in ${\rm Bimod}(A)$, with the comultiplication
given by the dual basis map and the counit given by the evaluation
map. Then ${\rm Hom}_{A}(D,A)_{A}={\rm End}(P)_{A}$. We call such
a representation {\em linear}.

\vspace{3mm}\paragraph{\textbf{Definition.}} Now we construct a
composition of convolution representations. Given two such
\begin{align}
B  \rightarrow {\rm Hom}_{A}(D,A)_{A},\ \ F \rightarrow {\rm
Hom}_{B}(E,B)_{B}.\nonumber
\end{align}
we define a third one.   It is a convolution representation given
as the composite of canonical ring homomorphisms
\begin{align}
E  \rightarrow {\rm Hom}_{B}(D,B)_{B}\rightarrow {\rm
Hom}_{B}(E,{\rm Hom}_{A}(D,A)_{A})_{B}={\rm
Hom}_{A}(E\otimes_{B^{o}\otimes B}D,A)_{A},\nonumber
\end{align}
where the structure of a left $(B^{o}\otimes B)$-module on $D$
comes from the $B$-bimodule structure defined in (\ref{bim}), and
the structure of a coalgebra in ${\rm Bimod}(A)$ \\ on
$E\otimes_{B^{o}\otimes B}D$ is defined as follows: the counits
$E\rightarrow B$ and $D\rightarrow A$ define a counit given as the
composite of canonical homomorphisms in ${\rm Bimod}(A)$
\begin{align}
E\otimes_{B^{o}\otimes B}D\rightarrow B\otimes_{B^{o}\otimes
B}D=D/[D,B]\rightarrow A,
\end{align}
the comultiplications
\begin{align}
D & \rightarrow D\otimes_{A}D,\ \ \ \ \  E\rightarrow E\otimes_{B}E,\nonumber\\
 d & \mapsto d_{(1)}\otimes d_{(2)},\ \ \ \   e\mapsto e_{(1)}\otimes e_{(2)},\nonumber
\end{align}
define a comultiplication
\begin{align}
E\otimes_{B^{o}\otimes B}D & \rightarrow (E\otimes_{B^{o}\otimes
B}D)\otimes_{A}(E\otimes_{B^{o}\otimes B}D),\nonumber\\
e\otimes d & \mapsto (e_{(1)}\otimes
d_{(1)})\otimes(e_{(2)}\otimes d_{(2)}).\nonumber
\end{align}
Identities have the following description
\begin{align}
{\rm id}_{A} : A \stackrel{=}{\rightarrow}{\rm Hom}_{A}(A\otimes
A,A)_{A}.\nonumber
\end{align}
In this way we obtain a category whose objects are (unital
associative) rings and morphisms are isoclasses of convolution
representations.
\begin{theorem}
There is a contravariant fully faithful embedding of the category
of rings with isoclasses of convolution representations as
morphisms into the category of noncommutative schemes.
\end{theorem}
{\em Proof.} To every associative ring $A$ we assign its monoidal
category $({\rm Bimod}(A), \otimes_{A}, A)$ of bimodules.  Given a
convolution representation
\begin{align}
B & \rightarrow {\rm Hom}_{A}(D,A)_{A},\label{rep}\\
\nonumber b & \mapsto (c\mapsto b(d)),
\end{align} for some coalgebra $D$ in ${\rm Bimod}(A)$ we have a canonical structure of a $B$-bimodule on $D$, commuting with the
original $A$-bimodule structure, defined as follows
\begin{align}
b\cdot d:=c_{(1)}b(d_{(2)}),\ \ \  d\cdot
b:=b(d_{(1)})d_{(2)}.\label{bim}
\end{align}
Regarding a $B$-bimodule structure as the same as a left or right
$B^{o}\otimes B$-module structure we obtain the following pair of
adjoint functors $f^{*}\dashv f_{*} $

\begin{align}
f^{*}: {\rm Bimod}(B) & \rightarrow {\rm Bimod}(A),\ \ \
 f^{*}N  =N\otimes_{B^{o}\otimes B}D,\\
 f_{*}: {\rm Bimod}(A) & \rightarrow {\rm Bimod}(B), \ \ \
f_{*}M ={\rm Hom}_{A}(D,M)_{A}.
\end{align}
For any $A$-bimodule $M$ below, we will denote by $N$ the
$B$-bimodule $f_{*}M ={\rm Hom}_{A}(D,M)_{A}$ and we will regard
elements of $N$ as maps $d\mapsto n(d)$. Then we have the
following natural transformation $f_{*}(-)\otimes_{B} f_{*}(-)
\rightarrow f_{*}(-\otimes_{A}-)$ of bi-functors ${\rm
Bimod}(A)\times {\rm Bimod}(A)\rightarrow {\rm Bimod}(B)$
\begin{align}
f_{*}M_{1}\otimes_{B} f_{*}M_{2} & \rightarrow
f_{*}(M_{1}\otimes_{A}M_{2}),\\
\nonumber (d\mapsto n_{1}(d))\otimes (d\mapsto n_{2}(d)) & \mapsto
(c\mapsto n_{1}(d_{(1)})\otimes n_{2}(d_{(2)})),
\end{align}
making $f_{*}$ a monoidal functor.

 Note that the
representation (\ref{rep}) gives rise automatically to a morphism
$B\rightarrow f_{*}A$ in ${\rm Bimod}(B)$. In this way a
convolution representation gives rise to a morphism of
noncommutative schemes.

Consider now a morphism of bimodule categories in the category of
noncommutative schemes given by a comonoidal functor $f^{*}: {\rm
Bimod}(B)\rightarrow {\rm Bimod}(A)$ and a homomorphism
$f^{*}B\rightarrow A$ in ${\rm Bimod}(A)$. Since $B\otimes B$ is a
coalgebra in ${\rm Bimod}(B)$, with the counit
\begin{align}B\otimes B & \rightarrow B\\
\nonumber b_{1}\otimes b_{2}& \rightarrow b_{1}b_{2}, \end{align}
and the comultiplication
\begin{align}B\otimes B & \rightarrow (B\otimes B)\otimes_{B}(B\otimes B)\\
\nonumber b_{1}\otimes b_{2} & \rightarrow (b_{1}\otimes
1)\otimes(1\otimes b_{2}), \end{align} $D:=f^{*}(B\otimes B)$ is a
coalgebra in ${\rm Bimod}(A)$ as well.  Let us consider the
composite of the following canonical maps, using the morphism
$f^{*}B\rightarrow A$ in ${\rm Bimod}(A)$
\begin{align}B & = {\rm Hom}_{B}(B\otimes B, B)_{B}\rightarrow
{\rm Hom}_{A}(f^{*}(B\otimes B),f^{*} (B))_{A}\\
\nonumber & \rightarrow {\rm Hom}_{A}(f^{*}(B\otimes
B),A)_{A}={\rm Hom}_{A}(D,A)_{A}.\end{align} Since $f^{*}$ is
comonoidal, the above composite is a homomomorphism of rings,
where ${\rm Hom}_{A}(D,A)_{A}$ is equipped with the canonical
convolution product.
 One can check easily that the functor $f^{*}$ is an image
 under our assignment of this convolution representation of $B$ over $A$. $\blacksquare$

\vspace{3mm}\paragraph{\textbf{Definition.}} We call the essential
image of the above embedding  \emph{category of affine
noncommutative schemes}. Let $A$ be an associative ring. Then we
define the object $X:={\rm Spec}(A)$ of $\mathfrak{Sch}$ as
follows: ${\rm Qcoh}(X):={\rm Bimod}(A)$,
$\otimes_{X}:=\otimes_{A}$, $\mathcal{O}_{X}:=A$. We call such an
affine noncommutative scheme \emph{(noncommutative) spectrum}. In
particular, we have a distinguished affine scheme $S:={\rm
Spec}(\mathbb{Z})$, where ${\rm Qcoh}(S):={\rm Ab}={\rm
Bimod}(\mathbb{Z})$, $\otimes_{S}:=\otimes=\otimes_{\mathbb{Z}}$,
$\mathcal{O}_{S}:=\mathbb{Z}$.

\vspace{3mm}

\subsubsection{Morita invariance of the spectrum} The following
fact is a corollary of the above structural theorem. Essentially,
it is a monoidal enhancement of the well known Morita invariance
of bimodule categories.

\begin{theorem}
A Morita equivalence between associative rings $A$ and $B$ induces
an isomorphism between affine noncommutative schemes ${\rm
Spec}(A)$ and ${\rm Spec}(B)$.
\end{theorem}

\emph{Proof.} A Morita equivalence of $A$ and $B$ can be described
as a representation which is an isomorphism of rings
\begin{align}
B\rightarrow {\rm Hom}_{A}(P^{*}\otimes P, A)_{A}={\rm End}(P)_{A}
,\label{morita}
\end{align}
where $P$ is a finitely generated projective generator in the
category $Mod$-$A$ of right $A$-modules. This defines a morphism
of affine schemes ${\rm Spec}(A)\rightarrow {\rm Spec}(B)$.
 By the Morita
theory $P^{*}$ is a finitely generated projective generator in
$Mod$-$B$ and the homomorphism of rings
\begin{align}
A\rightarrow{\rm Hom}_{A}(P^{*}\otimes P, A)_{A}={\rm
End}(P^{*})_{B},\label{in}
\end{align}
is an isomorphism. This defines a morphism ${\rm
Spec}(A)\leftarrow {\rm Spec}(B)$ in the opposite direction in the
same way as (\ref{morita}) defined $f$. It is inverse to $f$ since
by the Morita  theory
\begin{align}
P^{*}\otimes_{B}P\cong A\ \ ({\rm resp.}\ P\otimes_{A}P^{*}\cong
B) \end{align} in ${\rm Bimod}(A)$ (resp. in ${\rm Bimod}(B)$),
hence composites of these two morphisms are represented by
coalgebras
\begin{align}
(P\otimes P^{*})\otimes_{B^{o}\otimes B}(P^{*}\otimes P) &
=(P^{*}\otimes_{B} P)\otimes(P^{*}\otimes_{B} P)\cong A\otimes
A\nonumber\\
({\rm resp.}\ (P^{*}\otimes P)\otimes_{A^{o}\otimes A}(P\otimes
P^{*}) & =(P\otimes_{A} P^{*})\otimes(P\otimes_{A} P^{*})\cong
B\otimes B\nonumber)
\end{align} in  ${\rm Bimod}(A)$ (resp. in ${\rm Bimod}(B)$) isomorphic to coalgebras representing identity morphisms. $\blacksquare$

The problem of description of monoidal equivalences of categories
of bimodules was first considered by Takeuchi \cite{tak, tak'}
under the name of $\sqrt{Morita}$-equivalences. Since any strong
monoidal equivalence admits the left adjoint (which is necessarily
strong comonoidal) $\sqrt{Morita}$-equivalences are isomorphism in
the category of  affine noncommutative schemes, and all such
isomorphism are of that kind.

\subsection{Affine morphisms and spectra}

 \vspace{3mm}
\paragraph{\textbf{Definition.}} For every $A\in {\rm Alg}(X)$ and ${\mathcal F}_{1}$, ${\mathcal F}_{2}\in {\rm Bimod}_{X}(A)$ we
define the tensor product ${\mathcal F}_{1}\otimes_{A}{\mathcal
F}_{2}$, as usual, as the cokernel of the canonical pair of
morphisms
\begin{align}
{\mathcal F}_{1}\otimes_{X}A\otimes_{X}{\mathcal
F}_{2}\rightrightarrows {\mathcal F}_{1}\otimes_{X}{\mathcal
F}_{2}.
\end{align}

\vspace{3mm}
\paragraph{\textbf{Definition.}} We call a morphism $f: X\rightarrow
Y$ \emph{affine}  if the functor $f_{*}$ is faithful,  exact, and
the natural transformation of  bifunctors
\begin{align}f_{*}(-)\otimes_{f_{*}\mathcal{O}_{X}}f_{*}(-)\rightarrow
f_{*}(-\otimes_{X}-)\label{affiso}
\end{align}
is an isomorphism.

\vspace{3mm}
\paragraph{\textbf{Definition.}}
For every $A\in {\rm Alg}(X)$ such that the category $({\rm
Bimod}_{X}(A), \otimes_{A}, A)$ is abelian monoidal we define the
following noncommutative scheme
\begin{align}{\rm  Spec}_{X}(A):= & \ ({\rm Bimod}_{X}(A), \otimes_{A}, A)
\end{align}
and a canonical morphism  ${\rm Spec}_{X}(A)\rightarrow X$ whose
direct image functor is forgetting of the $A$-bimodule structure.

\vspace{3mm}
\paragraph{\textbf{Remark.}}
We have a canonical isomorphism of noncommutative schemes
\begin{align}
           X\rightarrow{\rm  Spec}_{X}(\mathcal{O}_{X}).
\end{align}

\begin{proposition}
Given an affine morphism $f: X\rightarrow Y$ the category \break
$({\rm Bimod}_{X}(f_{*}{\mathcal O}_{X}), \otimes_{f_{*}{\mathcal
O}_{X}}, f_{*}{\mathcal O}_{X})$ is abelian monoidal and we have
the following canonical decomposition
$$\begin{array}{ccc}
          & f         & \\
X & \overrightarrow{\ \ \ \ \ \ \ \ \ \ \ \ \ \ \ \ \ \ \ \ } & Y, \\
 \ \ \ \ \ \ \ \ \ \ \       \searrow  &                     & \nearrow\ \ \ \ \ \ \ \ \ \ \ \ \  \\
          & {\rm  Spec}_{Y}(f_{*}{\mathcal O}_{X})&
\end{array}$$
where the south-east arrow is an isomorphism.\label{decaff}
\end{proposition}

{\em Proof.} Since $f_{*}$ is monoidal it admits the following
canonical decomposition
$$\begin{array}{ccc}
          & f_{*}         & \\
{\rm Qcoh}(X) & \overrightarrow{\ \ \ \ \ \ \ \ \ \ \ \ \ \ \ \ \ \ \ \ \ \ \ \ \ \ \ } & {\rm Qcoh}(Y) \\
\ \ \ \ \ \ \ \ \ \ \ \ \ \ \ \ \ \ \ \ \ \ \ \ \       \searrow  &                     & \nearrow\ \ \ \ \ \ \ \ \ \ \ \ \ \ \ \ \ \ \ \ \ \ \  \\
          & {\rm Bimod}_{Y}(f_{*}{\mathcal O}_{X})&
\end{array}$$
where the right hand side arrow is the forgetting of the
$f_{*}{\mathcal O}_{X}$-bimodule structure.

Since $f_{*}$ is monoidal we have, for every two ${\mathcal
F}_{1}$, ${\mathcal F}_{2}\in {\rm Qcoh}(X)$, the following
decomposition
$$\begin{array}{ccc}
          f_{*}{\mathcal F}_{1}\otimes_{Y}f_{*}{\mathcal F}_{2} & \overrightarrow{\ \ \ \ \ \ \ \ \ \ \ \ \ \ \ \ \ \ \ \ } & f_{*}({\mathcal F}_{1}\otimes_{X}{\mathcal F}_{2}). \\
 \ \ \ \ \ \ \ \ \ \ \       \searrow  &                     & \nearrow\ \ \ \ \ \ \ \ \ \ \ \ \  \\
          & f_{*}{\mathcal F}_{1}\otimes_{f_{*}{\mathcal O}_{X}}f_{*}{\mathcal F}_{2}&
\end{array}$$
Since every equivalence of categories has the left adjoint (equal
to the inverse) and forgetting the bimodule structure has the left
adjoint, $f_{*}$ has the left adjoint as well. Since $f_{*}$ being
faithful exact induces by the Barr-Beck theorem an equivalence
${\rm Qcoh}(X) \rightarrow {\rm Bimod}_{Y}(f_{*}{\mathcal O}_{X})$
the latter category is abelian. $({\rm Bimod}_{Y}(f_{*}{\mathcal
O}_{X}), \otimes_{f_{*}{\mathcal O}_{X}},f_{*}{\mathcal O}_{X})$
is a monoidal category and the equivalence is strong monoidal by
(\ref{affiso}) applied to the latter decomposition. $\blacksquare$

\subsection{Affine schemes and affine morphisms} Exactly as in the
commutative algebraic geometry we have the following proposition.

\begin{proposition}
A noncommutative scheme is affine iff it admits an affine morphism
to ${\rm Spec}(\mathbb{Z})$.
\end{proposition}

\emph{Proof.} If $X={\rm Spec}(A)$ then the unique ring
homomorphism $\mathbb{Z}\rightarrow A$ defines a morphism
$X\rightarrow {\rm Spec}(\mathbb{Z})$. Since algebras over
$\mathbb{Z}$ are simply rings, then by Proposition \ref{decaff}
every affine morphism  $X\rightarrow {\rm Spec}(\mathbb{Z})$
defines an isomorphism
\begin{align}
X\stackrel{\cong}\rightarrow {\rm Spec}_{{\rm
Spec}(\mathbb{Z})}(f_{*}{\mathcal O}_{X})={\rm
Spec}(f_{*}{\mathcal O}_{X}).\ \blacksquare
\end{align}

 \vspace{3mm}
\paragraph{\textbf{Remark.}} Although every affine scheme ${\rm Spec}(A)$
admits a morphism ${\rm Spec}(A)\rightarrow{\rm Spec}(\mathbb{Z})$
corresponding to the unique ring homomorphism
$\mathbb{Z}\rightarrow A$, ${\rm Spec}(\mathbb{Z})$ is not a final
object in the category of noncommutative affine schemes because of
possible non-identical morphisms ${\rm
Spec}(\mathbb{Z})\rightarrow {\rm Spec}(\mathbb{Z})$ corresponding
to arbitrary coalgebras $D$ defined over $\mathbb{Z}$. This
resembles the situation in homotopy theory of topological
$G$-spaces, where the one-point space is not a final object and
there are spaces leaving under it, e.g. classifying spaces. In
this analogy coalgebras play the role of group actions. More
precise relation between coalgebras and group actions (or Hopf
algebra coactions) needs some regularity conditions discussed in
the next section.

\subsection{Flat covers and cospectra}

\vspace{3mm}
\paragraph{\textbf{Definition.}} For every $C\in {\rm Coalg}(X)$ and ${\mathcal F}_{1},{\mathcal
F}_{2}\in  {\rm Bicomod}_{X}(C)$ we define the cotensor product
${\mathcal F}_{1}\Box^{C}{\mathcal F}_{2}$, as usual, as the
kernel of the canonical pair of morphisms
\begin{align}
{\mathcal F}_{1}\otimes_{X}{\mathcal
F}_{2}\rightrightarrows{\mathcal
F}_{1}\otimes_{X}C\otimes_{X}{\mathcal F}_{2}.
\end{align}

\vspace{3mm}
\paragraph{\textbf{Definition.}} We call a morphism $f: X\rightarrow
Y$ \emph{flat }  if the functor $f^{*}$ is exact  and the natural
transformation of  bifunctors
\begin{align}f^{*}(-)\Box^{f^{*}\mathcal{O}_{Y}}f^{*}(-)\leftarrow
f^{*}(-\otimes_{Y}-)\label{flatiso}
\end{align}
is an isomorphism, and {\em cover} if $f^{*}$ is faithful.

\vspace{3mm}
\paragraph{\textbf{Definition.}}
For every $C\in {\rm Coalg}(X)$ such that the category $({\rm
Bicomod}_{X}(C), \Box^{C}, C)$ is abelian monoidal we define the
following noncommutative scheme
\begin{align}
           {\rm  Cospec}_{X}(C):=  \ ({\rm  Bicomod}_{X}(C), \Box^{C}, C).
\end{align}
and a canonical morphism  $X\rightarrow {\rm Cospec}_{X}(C)$,
whose inverse image functor is forgetting of the $C$-bicomodule
structure.

\vspace{3mm}
\paragraph{\textbf{Remark.}}
We have a canonical isomorphism of noncommutative schemes
\begin{align}
           {\rm  Cospec}_{X}(\mathcal{O}_{X})\rightarrow X.
\end{align}

\begin{proposition}
Given a flat cover $f: X\rightarrow Y$ the category \break $({\rm
Bicomod}_{X}(f^{*}{\mathcal O}_{Y}), \Box^{f^{*}{\mathcal O}_{Y}},
f^{*}{\mathcal O}_{Y})$ is abelian monoidal and we have the
following canonical decomposition
$$\begin{array}{ccc}
          & f         & \\
X & \overrightarrow{\ \ \ \ \ \ \ \ \ \ \ \ \ \ \ \ \ \ \ \ } & Y, \\
 \ \ \ \ \ \ \ \ \ \ \       \searrow  &                     & \nearrow\ \ \ \ \ \ \ \ \ \ \ \ \  \\
          & {\rm  Cospec}_{X}(f^{*}{\mathcal O}_{Y})&
\end{array}$$
where the north-east arrow is an isomorphism.
\end{proposition}

{\em Proof.} Since $f^{*}$ is comonoidal it admits the following
canonical decomposition
$$\begin{array}{ccc}
          & f^{*}         & \\
{\rm Qcoh}(X) & \overleftarrow{\ \ \ \ \ \ \ \ \ \ \ \ \ \ \ \ \ \ \ \ \ \ \ \ \ \ \ } & {\rm Qcoh}(Y) \\
\ \ \ \ \ \ \ \ \ \ \ \ \ \ \ \ \ \ \ \ \ \ \ \ \       \nwarrow  &                     & \swarrow\ \ \ \ \ \ \ \ \ \ \ \ \ \ \ \ \ \ \ \ \ \ \  \\
          & {\rm Bicomod}_{X}(f^{*}{\mathcal O}_{Y})&
\end{array}$$
where the north-west arrow is  forgetting of the $f^{*}{\mathcal
O}_{Y}$-bicomodule structure.

Since $f^{*}$ is comonoidal we have, for every two ${\mathcal
G}_{1}$, ${\mathcal G}_{2}\in {\rm Qcoh}(Y)$, the following
decomposition
$$\begin{array}{ccc}
          f^{*}{\mathcal G}_{1}\otimes_{X}f^{*}{\mathcal G}_{2} & \overleftarrow{\ \ \ \ \ \ \ \ \ \ \ \ \ \ \ \ \ \ \ \ } & f^{*}({\mathcal G}_{1}\otimes_{Y}{\mathcal G}_{2}). \\
 \ \ \ \ \ \ \ \ \ \ \       \nwarrow  &                     & \swarrow\ \ \ \ \ \ \ \ \ \ \ \ \  \\
          & f^{*}{\mathcal G}_{1}\Box^{f^{*}{\mathcal O}_{Y}}f^{*}{\mathcal G}_{2}&
\end{array}$$
Since every equivalence of categories has the right adjoint (equal
to the inverse) and forgetting the bicomodule structure has the
right adjoint, $f^{*}$ has the right adjoint as well. Since
$f^{*}$ being faithful exact   induces by the Barr-Beck theorem an
equivalence ${\rm Qcoh}(Y) \leftarrow {\rm
Bicomod}_{X}(f^{*}{\mathcal O}_{Y})$  the latter category  is
abelian. $({\rm Bicomod}_{X}(f^{*}{\mathcal O}_{Y}),
\Box^{f_{*}{\mathcal O}_{Y}},f^{*}{\mathcal O}_{Y})$ is a monoidal
category and the equivalence is strong comonoidal by
(\ref{flatiso}) applied to the latter decomposition.
$\blacksquare$

\subsection{Sections and the sheaf condition} Consider a commutative
diagram
$$\begin{array}{ccccc}
         & \pi         &     \\
U     & \overrightarrow{\ \ \ \ \ \ \ } &  X, \\
 \ \ \ \ \ \ \ \ \ \ \ \ \        \alpha  \searrow  &        &    \swarrow \beta  \ \ \ \ \ \ \ \ \ \ \ \ \  \\
         & S &
\end{array}$$
in the category of noncommutative schemes over $S:={\rm
Spec}(\mathbb{Z})$ with $\pi$  flat . Let ${\mathcal{F}}\in {\rm
Qcoh}(X)$.

\vspace{3mm}
\paragraph{\textbf{Definition.}} We define an abelian {\em group of
sections of ${\mathcal{F}}$ over $U$}  as follows
\begin{align}
{\mathcal{F}}(U):={\rm Hom}_{{\rm
Cospec}_{X}(\pi^{*}\mathcal{O}_{X})}(\alpha^{*}\mathcal{O}_{S},
\pi^{*}{\mathcal{F}}).
\end{align}

\begin{proposition}
${\mathcal O}_{X}(-)$ is a contravariant  functor from the
category of $S$-schemes flat over $X$ to the category of rings.
For any ${\mathcal{F}}\in {\rm Qcoh}(X)$ sections
${\mathcal{F}}(-)$ form an ${\mathcal O}_{X}(-)$-bimodule and the
assignment ${\mathcal{F}}\rightsquigarrow{\mathcal{F}}(-)$ is a
monoidal functor ${\rm Qcoh}(X)\rightarrow {\rm Bimod}({\mathcal
O}_{X}(-))$.
\end{proposition}

{\em Proof.}  Since $\beta^{*}\mathcal{O}_{S}\in {\rm Coalg}(X)$
and $\pi$ is flat the inverse image
$\alpha^{*}\mathcal{O}_{S}=\pi^{*}\beta^{*}\mathcal{O}_{S}$ is a
coalgebra in the monoidal category $({\rm
Bicomod}(\pi^{*}{\mathcal O}_{X}), \Box^{\pi^{*}{\mathcal O}_{X}},
\pi^{*}{\mathcal O}_{X})$. Then the following canonical composite
\begin{align}
{\rm Hom}^{ \pi^{*}{\mathcal
O}_{X}}_{U}(\alpha^{*}\mathcal{O}_{S}, \pi^{*}{\mathcal O}_{X})^{
\pi^{*}{\mathcal O}_{X}}& \otimes{\rm
Hom}^{ \pi^{*}{\mathcal O}_{X}}_{U}(\alpha^{*}\mathcal{O}_{S}, \pi^{*}{\mathcal O}_{X})^{ \pi^{*}{\mathcal O}_{X}}\label{ringstr}\\
& \downarrow \nonumber\\
{\rm Hom}^{ \pi^{*}{\mathcal O}_{X}}_{U}(\alpha^{*}\mathcal{O}_{S}\Box^{ \pi^{*}{\mathcal O}_{X}}\alpha^{*}\mathcal{O}_{S},& \pi^{*}{\mathcal O}_{X}\Box^{\pi^{*}{\mathcal O}_{X}}\pi^{*}{\mathcal O}_{X})^{\pi^{*}{\mathcal O}_{X}}\nonumber\\ & \downarrow\nonumber\\
{\rm Hom}^{\pi^{*}{\mathcal O}_{X}}_{U}(\alpha^{*}\mathcal{O}_{S}&
,\pi^{*}{\mathcal O}_{X})^{\pi^{*}{\mathcal O}_{X}}.\nonumber
\end{align}
defines a ring structure on sections ${\mathcal O}_{X}(U)={\rm
Hom}_{{\rm
Cospec}_{X}(\pi^{*}\mathcal{O}_{X})}(\alpha^{*}\mathcal{O}_{S},
\pi^{*}{\mathcal O}_{X})={\rm Hom}_{U}^{\pi^{*}{\mathcal
O}_{X}}(\alpha^{*}\mathcal{O}_{S},\pi^{*}{\mathcal
O}_{X})^{\pi^{*}{\mathcal O}_{X}}$ with the unit being the image
of the identity in ${\rm Hom}_{U}^{\pi^{*}{\mathcal
O}_{X}}(\pi^{*}{\mathcal O}_{X},\pi^{*}{\mathcal
O}_{X})^{\pi^{*}{\mathcal O}_{X}}$ under the canonical map induced
by the morphism
$\alpha^{*}\mathcal{O}_{S}=\pi^{*}\beta^{*}\mathcal{O}_{S}\rightarrow
\pi^{*}{\mathcal O}_{X}$. Replacing in (\ref{ringstr}) one copy of
${\mathcal O}_{X}$ in the covariant argument of Hom by
${\mathcal{F}}$ we obtain an ${\mathcal O}_{X}(U)$-bimodule
structure on ${\mathcal{F}}(U)$. Similarly, replacing in
(\ref{ringstr}) two copies of ${\mathcal O}_{X}$ by
$\mathcal{F}_{1}$, $\mathcal{F}_{2}$ and inverting the isomorphism
\begin{align}
\pi^{*}(\mathcal{F}_{1})\Box^{\pi^{*}\mathcal{O}_{X}}\pi^{*}(\mathcal{F}_{2})\leftarrow
\pi^{*}(\mathcal{F}_{1}\otimes_{X}\mathcal{F}_{2})
\end{align}
we obtain
\begin{align}
\mathcal{F}_{1}(U)\otimes\mathcal{F}_{2}(U)\rightarrow
(\mathcal{F}_{1}\otimes_{X}\mathcal{F}_{2})(U).
\end{align}
One can check that the above map factorizes canonically through
\begin{align}
\mathcal{F}_{1}(U)\otimes_{\mathcal{O}_{X}(U)}\mathcal{F}_{2}(U)\rightarrow
(\mathcal{F}_{1}\otimes_{X}\mathcal{F}_{2})(U).
\end{align}
It is clear that all these constructions are functorial in $U$. $
\blacksquare$

\begin{proposition} If $\pi$ is a flat cover then
${\mathcal{F}}(U)={\mathcal{F}}(X)$. In particular, the functor of
global sections is independent of the choice of a flat cover.
\end{proposition}

{\em Proof.} By Proposition 2  we have
\begin{align}
{\mathcal{F}}(U)& ={\rm Hom}_{{\rm
Cospec}_{X}(\pi^{*}\mathcal{O}_{X})}(\pi^{*}\beta^{*}\mathcal{O}_{S},
\pi^{*}{\mathcal{F}}) ={\rm Hom}_{X}(\beta^{*}\mathcal{O}_{S},
{\mathcal{F}})={\mathcal{F}}(X).\ \blacksquare\nonumber
\end{align}

\begin{proposition}
${\mathcal{F}}(X)=\beta_{*}{\mathcal{F}}$.
\end{proposition}

{\em Proof.} By the following isomorphism of functors   ${\rm
Hom}_{S}(\mathcal{O}_{S}, -)=id_{{\rm Ab}}$, for $S={\rm
Spec}(\mathbb{Z})$, we have
\begin{align}
{\mathcal{F}}(X) ={\rm Hom}_{X}(\beta^{*}\mathcal{O}_{S},
{\mathcal{F}}) ={\rm Hom}_{S}(\mathcal{O}_{S},
\beta_{*}{\mathcal{F}})=\beta_{*}{\mathcal{F}}.\
\blacksquare\nonumber
\end{align}

\subsection{Flat covers and Galois extensions}
Let $A$ be a ring and $D\in {\rm Coalg}(A)$.

\vspace{3mm}
\paragraph{\textbf{Definition.}} Let $R$ be a ring. A morphism $\varphi: M_{0}\rightarrow
M_{1}$ in ${\rm Bimod}(R)$ is called \emph{pure} if for all $U\in
{\rm Mod}_{R}$, $V\in _{R}{\rm Mod}$ the induced sequence
\begin{align}0\rightarrow U\otimes_{R}{\rm ker}(\varphi)\otimes_{R}V\rightarrow U\otimes_{R}M_{0}\otimes_{R}V\rightarrow
U\otimes_{R}M_{1}\otimes_{R}V\end{align} is exact.

The following lemma is a bimodule version of the purity criterion
of compatibility  of the tensor and cotensor products
\cite{tak'',schn}.

\begin{lemma} Assume that there is given $C\in{\rm Coalg}(A)$ and a left
$R$-module right $C$-comodule $M_{1}$ and a right  $R$-module left
$C$-comodule $M_{2}$. Then

1) for every $U\in {\rm Mod}_{R}$ and every $V\in _{R}{\rm Mod}$
there exists a canonical, natural in $U$ and $V$, homomorphism of
abelian groups
\begin{align}
U\otimes_{R}(M_{1}\square^{C}_{A}M_{2})\otimes_{R}V\rightarrow
(U\otimes_{R}M_{1})\square^{C}_{A}(M_{2}\otimes_{R}V)\label{ab}
\end{align}

2) the canonical morphism in  ${\rm Bimod}(R)$ defining the
cotensor product
\begin{align}
M_{1}\otimes_{A}M_{2}\rightarrow M_{1}\otimes_{A}C\otimes_{A}M_{2}
\end{align}
is pure iff all homomorphisms (\ref{ab}) are
isomorphisms.\label{lemma1}
\end{lemma}

\emph{Proof.} We have a canonical exact sequence in ${\rm
Bimod}(R)$
\begin{align}0\rightarrow M_{1}\square^{C}_{A}M_{2}\rightarrow M_{1}\otimes_{A}M_{2}\rightarrow M_{1}\otimes_{A}C\otimes_{A}M_{2}.
\end{align}

Tensoring it from both sides by $U$ and $V$ we obtain a complex
fitting into the following diagram with the exact bottom row and
the left hand side arrow in the bottom row injective
$$\begin{array}{ccccc}
U\!\otimes_{R}\!(M_{1}\square^{C}_{A}M_{2})\!\otimes_{R}\! V &
\!\!\! \rightarrow U\!\otimes_{R}\!(M_{1}\!\otimes_{A}\!
M_{2})\!\otimes_{R}\! V &\!\!\! \rightarrow
U\!\otimes_{R}\!(M_{1}\!\otimes_{A}\! C\!\otimes_{A}\!
M_{2})\!\otimes_{R}\! V\\
\downarrow &  \ \ \ \downarrow \cong & \ \ \ \ \ \downarrow \cong\\
(U\!\otimes_{R}\! M_{1})\square_{A}^{C}(\! M_{2}\!\otimes_{R}\! V)
& \!\!\! \rightarrowtail (U\!\otimes_{R}\! M_{1})\!\otimes_{A}\!
(M_{2}\!\otimes_{R}\! V) &\!\!\! \rightarrow (U\!\otimes_{R}\!
M_{1})\!\otimes_{A}\! C\!\otimes_{A}\! (M_{2}\!\otimes_{R}\! V)
\end{array}$$
This defines the vertical left hand side arrow and implies that
the upper row is exact  and the left hand side arrow in the upper
row is injective iff the vertical left hand side arrow is an
isomorphism.\ $\blacksquare$

\vspace{3mm}
\paragraph{\textbf{Definition.}} If there is given a morphism  $D\rightarrow C$ in ${\rm
Coalg}(A)$ such that $\Box^{C}_{A}$ is associative (with the unit
$C$) then $D$ becomes a coalgebra in the monoidal category $({\rm
Bicomod}_{A}(C), \Box^{C}_{A}, C )$ and the canonical composite
\begin{align}
{\rm Hom}^{C}_{A}(D,C)^{C}_{A}\otimes{\rm
Hom}^{C}_{A}(D,C)^{C}_{A}\rightarrow {\rm
Hom}^{C}_{A}(D\Box^{C}_{A}D,C\Box^{C}_{A}C)^{C}_{A}\rightarrow{\rm
Hom}^{C}_{A}(D,C)^{C}_{A},
\end{align}
where the right hand arrow uses the comultiplication $D\rightarrow
D\Box^{C}_{A}D$ and the isomorphism
$C\stackrel{\cong}{\rightarrow} C\Box^{C}_{A}C$, defines a subring
structure on ${\rm Hom}^{C}_{A}(D,C)^{C}_{A}\subset {\rm
Hom}_{A}(D,A)_{A}$, which we call {\em subring of invariants}.

\vspace{3mm}
\paragraph{\textbf{Definition.}} Assume there is given a
representation $B\rightarrow {\rm Hom}_{A}(D,A)_{A}$ factorizing
through the subring of invariants  ${\rm
Hom}^{C}_{A}(D,C)^{C}_{A}$. It is an exercise in Sweedler's
notation to prove using (\ref{bim}) that

\begin{itemize}
\item $D/[D,B]$ inherits a structure of a coalgebra in ${\rm Bimod}(A)$ from $D$,
\item the homomorphism $D\rightarrow C$ in ${\rm Coalg}(A)$
factorizes in ${\rm Coalg}(A)$ canonically through
\begin{align}
D/[D,B]\rightarrow C,
\end{align}
\item there are
structures of $(B^{o}\otimes B)$-bimodules on $B\otimes D$ and
$D\square^{C}_{A}D$ defined as follows
\begin{align}
(b'^{o}\otimes b'')\cdot(b\otimes d) & :=bb'\otimes b''\cdot d,\\
(b\otimes d)\cdot (b'^{o}\otimes b'') & :=b'b\otimes d\cdot b'',\nonumber\\
(b'^{o}\otimes b'')\cdot(d'\otimes d'') & :=d'\otimes b''\cdot
d''\cdot b' ,\\
(d'\otimes d'')\cdot(b'^{o}\otimes b'') & :=b'\cdot d'\cdot
b''\otimes d''.\nonumber
\end{align}
and there is a canonical morphism in ${\rm Bimod}(B^{o}\otimes B)$
\begin{align} B\otimes D\rightarrow D\square^{C}_{A}D,\nonumber\end{align}
defined as follows
\begin{align} b\otimes d\mapsto d_{(1)}b(d_{(2)})\otimes d_{(3)}=d_{(1)}\otimes b(d_{(2)})d_{(3)}.\end{align}
\item for every $M\in {\rm Bimod}(A)$ there is a canonical
morphism in ${\rm Bimod}(A)$ natural in $M$
\begin{align}
{\rm Hom}_{A}(D,M)_{A}\otimes_{B^{o}\otimes B }D\rightarrow
D/[D,B]\otimes_{A}M\otimes_{A}D/[D,B],\nonumber
\end{align}
defined as follows
\begin{align}
\mu\otimes d\mapsto [d_{(1)}]\otimes \mu(d_{(2)})\otimes
[d_{(3)}],
\end{align}
where $[d]:=d+[D,B]$.
\end{itemize}

\vspace{3mm}
\paragraph{\textbf{Definition.}}  We call a representation $B\rightarrow {\rm
Hom}_{A}(D,A)_{A}$ {\em noncommutative Galois ring extension} if
there is given a morphism  $D\rightarrow C$ in ${\rm Coalg}(A)$
such that
\begin{itemize}
\item ({\em regularity}) the category $({\rm Bicomod}_{A}(C), \Box^{C}_{A},C)$ is abelian monoidal,
\item ({\em purity}) the canonical morphism in ${\rm Bimod}(B^{o}\otimes B)$
\begin{align} D\otimes_{A}D\rightarrow D\otimes_{A}C\otimes_{A}D\nonumber\end{align}
defining the
cotensor product is pure,
\item  ({\em invariants}) $B$ is
 mapped isomorphically onto the {\em subring of invariants}
 \begin{align}{\rm
Hom}^{C}_{A}(D,C)^{C}_{A}\subset {\rm
Hom}_{A}(D,A)_{A},\nonumber\end{align} and the canonical morphism
in ${\rm Bimod}(B^{o}\otimes B)$
\begin{align} B\otimes D\rightarrow D\square^{C}_{A}D,\nonumber\end{align}
is an isomorphism,
\item ({\em faithful flatness}) $D$ is faithfully flat as a left $B^{o}\otimes B$-module,
\item ({\em freeness}) the {\em canonical morphism}
\begin{align} D/[D,B]\rightarrow C\nonumber\end{align} in ${\rm Coalg}(A)$ is an
isomorphism,
\item ({\em comonad}) the natural transformation of functors
\begin{align}
{\rm Hom}_{A}(D,-)_{A}\otimes_{B^{o}\otimes B }D\rightarrow
C\otimes_{A}-\otimes_{A}C\nonumber
\end{align}
is an isomorphism.
\end{itemize}

\subsubsection{Comparison with Galois comodules.} Let $A$
be a ring, $C\in{\rm Coalg}(A)$ and $P$ be a finitely generated
projective right $A$-module. Assume now in addition that $P$ is  a
right comodule over $C$. Define $B:={\rm End}(P)^{C}_{A}$. Then
$P$ is called \emph{Galois comodule} if the natural transformation
of functors ${\rm Mod}_{A}\rightarrow{\rm Comod}^{C}_{A}$
\begin{align}
{\rm Hom}(P,-)_{A}\otimes_{B}P\rightarrow -\otimes_{A}C
\end{align}
is an isomorphism. For the complicated history of this simple
definition we refer the reader to \cite{wis'} (where the
assumption of beeing finitely generated projective is dropped).
Under some assumptions (\cite{wis'} Theorem 5.7. Equivalences.)
the functor
\begin{align}
{\rm Hom}(P,-)^{C}_{A}: {\rm Comod}^{C}_{A}\rightarrow {\rm
Mod}_{B}
\end{align}
is an equivalence with the inverse $-\otimes_{B}P$.

Galois comodules and their predecessors (Galois corings, Galois
extensions, coalgebra-Galois extensions, Hopf-Galois extensions,
\cite{brz,brz',brz'',brzhaj,brzwis,elkgom,schsch,wis}) found an
interesting interpretation in terms of descent theory, theory of
invariants, associated vector bundles etc. in the realm of
noncommutative geometry modeled on one-sided (co)modules
\cite{brz''}. Starting from slightly different assumptions, we
realize a similar program in our noncommutative geometry modeled
on monoidal categories. The main difference consists in the role
of the monoidal structure. In the first approach associativity of
the cotensor product fails in general (\cite{gru}) and can be
achieved only for the price of painful and difficult to examine
assumptions (\cite{gom}). In our approach cotensor products
arising from our  geometric flat covers are associative almost by
definition, so we put this associativity as a necessary regularity
condition in our definition of a noncommutative Galois ring
extension. In the purily commutative case of schemes with trivial
group scheme action, viewed as symmetric monoidal categories, the
inverse image of the monoidal unit $f^{*}{\mathcal
O}_{Y}={\mathcal O}_{X}$ is a monoidal unit again, so the
associativity condition is void. In the sequel of this paper we
will present the derived version of our monoidal noncommutative
geometry, where this restrongion disappears.  In general, it is
not easy to analyse relations between different assumptions in
one-sided and monoidal noncommutative geometries. The following
proposition can be regarded as a transition of the border line
between these two approaches.

\begin{proposition} Let $C\in{\rm Coalg}(A)$ be a coalgebra over a
ring $A$ such that $({\rm Bicomod}_{A}(C), \square^{C}_{A}, C)$ is
an abelian monoidal  category. Let $P$ be a right $C$-comodule
which is finitely generated and projective as a right $A$-module
and define a subring $B:={\rm End}(P)^{C}_{A}\subset {\rm
End}(P)_{A}$. Assume the following conditions are fulfilled:

1) the functor $P^{*}\otimes_{B}-\otimes_{B}P: {\rm
Bimod}(B)\rightarrow {\rm Bimod}(A)$ is faithful exact,

2) the first differential in the Amitsur complex
\begin{align}
{\rm End}(P)_{A}\rightarrow {\rm End}(P)_{A}\otimes_{B}{\rm
End}(P)_{A},\label{ami}
\end{align}
\begin{align}
e\mapsto 1\otimes e-e\otimes 1,\nonumber
\end{align}
is pure in ${\rm Bimod}(B)$,

3) the canonical morphism in ${\rm Coalg}(A)$
\begin{align}P^{*}\otimes_{B}P\rightarrow C,\label{fre}\end{align}
\begin{align}p^{*}\otimes p\mapsto p^{*}(p_{(0)})\cdot p_{(1)},\nonumber\end{align}
is an isomorphism.

Take $D:=P^{*}\otimes P$ and the canonical morphism $D\rightarrow
C$ in ${\rm Coalg}(A)$ defined as in (\ref{fre}). Then there is a
canonical ring isomorphism ${\rm End}(P)_{A}\cong {\rm
Hom}_{A}(D,A)_{A}$ such that $(B\rightarrow {\rm
Hom}_{A}(D,A)_{A}, D\rightarrow C)$ is a noncommutative Galois
ring extension.\label{lingal}
\end{proposition}

{\em Proof}. Since $P$ is finitely generated projective there is a
canonical isomorphism of rings
\begin{align}{\rm End}(P)_{A}\stackrel{\cong}{\rightarrow}{\rm Hom}_{A}(P^{*}\otimes P,A)_{A}\label{end},\end{align}
\begin{align}e\mapsto (p^{*}\otimes p\mapsto p^{*}(e(p))),\nonumber\end{align}
with the inverse defined by means of the dual basis
$((p_{i})_{i\in I}, (p^{*}_{i})_{i\in I})$,
$p^{*}_{i}(p_{j})=\delta_{ij}$, as follows
\begin{align}(p\mapsto \sum_{i\in I}p_{i}\cdot h(p^{*}_{i}\otimes p))\leftarrow\!\shortmid h,\end{align}

For a right $C$-comodule $P$, with the comultiplication
\begin{align}P\rightarrow P\otimes_{A}C, \end{align}
\begin{align}p\mapsto p_{(0)}\otimes p_{(1)},\nonumber \end{align}
which is finitely generated and projective as a right $A$-module,
the canonical left $C$-comodule structure on the dual $P^{*}:={\rm
Hom}(P,A)_{A}$ can be defined  by means of the dual basis
 as follows
\begin{align}P^{*}\rightarrow C\otimes_{A}P^{*}, \end{align}
\begin{align}p^{*}\mapsto p^{*}_{(-1)}\otimes p^{*}_{(0)}:=\sum_{i\in I}p^{*}(p_{i(0)})\cdot p_{i(1)}\otimes p^{*}_{i}. \end{align}
Then the canonical isomorphism (\ref{end}) induces an isomorphism
\begin{align}{\rm End}(P)^{C}_{A}\stackrel{\cong}{\rightarrow}{\rm Hom}^{C}_{A}(P^{*}\otimes P,A)_{A}^{C}\label{invend},\end{align}
which, by the definition of $B$, proves the first part of the
\emph{invariants} condition.

By the canonical isomorphism
\begin{align}P\otimes_{A}P^{*}\rightarrow {\rm End}(P)_{A}\label{dualb},\end{align}
\begin{align}p\otimes p^{*}\mapsto (p'\mapsto p\cdot p^{*}(p')),\nonumber\end{align}
with the inverse (dual basis)
\begin{align}\sum_{i\in I}e(p_{i})\otimes p^{*}_{i}\leftarrow\!\shortmid e,\end{align}
and the isomorphism (\ref{fre})
 the first
differential in the Amitsur complex is isomorphic to the map
defining the cotensor product $P\square^{C}_{A}P^{*}$
\begin{align}P\otimes_{A}P^{*}\rightarrow
P\otimes_{A}C\otimes_{A}P^{*},\label{purP}\end{align}
\begin{align}p\otimes p^{*}\mapsto
p_{(0)}\otimes p_{(1)}\otimes p^{*}-p\otimes p^{*}_{(-1)}\otimes
p^{*}_{(0)}.\nonumber\end{align}

By purity of (\ref{purP}) and Lemma \ref{lemma1}, the homomorphism
\begin{align} B\otimes D\rightarrow D\square^{C}_{A}D,\label{fir}\end{align}
can be rewritten as follows
\begin{align} B\otimes P^{*}\otimes P \rightarrow (P^{*}\otimes P)\square^{C}_{A}(P^{*}\otimes P)
                                      \cong P^{*}\otimes (P\square^{C}_{A}P^{*})\otimes P.\label{compat}\end{align}
The isomorphism (\ref{dualb}) induces an isomorphism
\begin{align}P\square^{C}_{A}P^{*}\rightarrow {\rm
End}(P)^{C}_{A},\label{cot}\end{align} hence (\ref{compat}) can be
rewritten as
\begin{align} B\otimes P^{*}\otimes P \rightarrow P^{*}\otimes {\rm
End}(P)^{C}_{A}\otimes P \cong {\rm End}(P)^{C}_{A}\otimes
P^{*}\otimes P\end{align} which is induced by the homomorphism
$B\rightarrow {\rm End}(P)^{C}_{A}$. If the latter homomorphism is
an isomorphism, (\ref{fir}) is an isomorphism as well, which
proves the second part of the \emph{invariants} condition.

Again by purity of (\ref{purP}) for all $U,V\in {\rm Bimod}(B)$
the canonical sequence
$$\begin{array}{c}
0\\
\downarrow\\
(P^{*}\otimes_{B}U)\otimes_{B}(P\square^{C}_{A}P^{*})\otimes_{B}(V\otimes_{B}P)\\
\downarrow\\
(P^{*}\otimes_{B}U)\otimes_{B}(P\otimes_{A}P^{*})\otimes_{B}(V\otimes_{B}P)\\
\downarrow\\
(P^{*}\otimes_{B}U)\otimes_{B}(P\otimes_{A}C\otimes_{A}P^{*})\otimes_{B}(V\otimes_{B}P),
\end{array}$$
is exact. It can be rewritten as
$$\begin{array}{c}
0\\
\downarrow\\
U\otimes_{B^{o}\otimes B}(P^{*}\otimes(P\square^{C}_{A}P^{*})\otimes P)\otimes_{B^{o}\otimes B}V\\
\downarrow\\
U\otimes_{B^{o}\otimes B}(P^{*}\otimes(P\otimes_{A}P^{*})\otimes P)\otimes_{B^{o}\otimes B}V\\
\downarrow\\
U\otimes_{B^{o}\otimes
B}(P^{*}\otimes(P\otimes_{A}C\otimes_{A}P^{*})\otimes
P)\otimes_{B^{o}\otimes B}V,
\end{array}$$
and further, using purity of (\ref{purP}) together with Lemma
\ref{lemma1}, and finally the definition of $D:=P^{*}\otimes P$,
it can be rewritten as
$$\begin{array}{c}
0\\
\downarrow\\
U\otimes_{B^{o}\otimes B}(D\square^{C}_{A}D)\otimes_{B^{o}\otimes B}V\\
\downarrow\\
U\otimes_{B^{o}\otimes B}(D\otimes_{A}D)\otimes_{B^{o}\otimes B}V\\
\downarrow\\
U\otimes_{B^{o}\otimes
B}(D\otimes_{A}C\otimes_{A}D)\otimes_{B^{o}\otimes B}V,
\end{array}$$
which proves purity of
\begin{align}
D\otimes_{A}D\rightarrow
D\otimes_{A}C\otimes_{A}D,\label{comult}\end{align} i.e. the
\emph{purity} condition.

The identification
\begin{align}(P^{*}\otimes P)/[P^{*}\otimes P,B]=P^{*}\otimes_{B} P
\end{align}
and the assumption 3) prove the \emph{freeness} condition.

By duality between right and left finitely generated projective
modules we have
\begin{align}
{\rm Hom}_{A}(P^{*}\otimes P, -)_{A}\otimes_{B^{o}\otimes B
}(P^{*}\otimes P)=(P^{*}\otimes_{B} P)
\otimes_{A}-\otimes_{A}(P^{*}\otimes_{B} P), \end{align} which
proves the {\em comonad} condition.\ $\blacksquare$


\vspace{3mm}
\paragraph{\textbf{Remark.}} If $K$ is a commutative ring, $H$ a Hopf algebra over $K$, $A$ a right comodule
algebra over $H$, then $C := A\otimes_{K}H \in {\rm Coalg}(A)$, $P
:= A \in {\rm Comod}^{C}_{A}$, $B := {\rm
Hom}^{C}_{A}(D,C)^{C}_{A} = A^{co\ H}$, and the canonical map
\begin{align}A\otimes_{B}A \rightarrow A\otimes_{K} H\end{align}
form a Hopf-Galois context. One can also take instead of a Hopf
algebra an arbitrary (symmetric) coalgebra over $K$ \cite{brz,
brzhaj} and form a coalgebra-Galois context.

If $S$ is a commutative affine scheme, $G$ a commutative group
scheme flat affine over $S$, acting freely
$U\times_{S}G\rightarrow U$ on a commutative scheme $U$ flat
affine over $S$ with a good quotient $X = U//G$, then $K :=
{\mathcal O}(S)$, $H := {\mathcal O}(G)$, $A := {\mathcal O}(U)$,
and $B ={\mathcal O}(X)={\mathcal O}(U//G)=A^{G}$ form a
Hopf-Galois context.

In particular, if $K$ is a field, $G$ a finite group of
automorphisms of a finite field extension $K\subset A$ then the
Hopf-Galois context with  $H := {\rm Map}(G,K)$ is equivalent to
the classical $G$-Galois field extension context
\begin{align}K\stackrel{\cong}{\rightarrow}A^{G}\Leftrightarrow A\otimes_{K} A \stackrel{\cong}{\rightarrow}{\rm
Map}(G,A).\end{align}

\vspace{3mm}
\paragraph{\textbf{Remark.}} The purity assumption in Proposition
\ref{lingal} is satisfied if the Amitsur complex admits a homotopy
contracting it to $B$, i.e. if $B$ is a direct summand in the
$B$-bimodule ${\rm End}(P)_{A}$. By the \emph{invariants}
condition the latter property is equivalent to the existence of a
generalized Reynolds operator, i.e. to generalized linear
reductivity of the coaction of $C$.

\subsubsection{Flat covers and Galois extensions.} The next theorem is the main result of this paper. In
particular, it means that various global constructions in
commutative geometry related to group actions and gluing
(invariants, descent theory) can be viewed as natural
constructions in noncommutative geometry.

\begin{theorem} There is one-to-one correspondence
between flat covers in the category of noncommutative affine
schemes and noncommutative Galois ring extensions.
\end{theorem}
{\em Proof.} ({\em Flat covers $\rightsquigarrow$ Galois
extensions}) Consider a flat cover $f:{\rm Spec}(A)\rightarrow
{\rm Spec}(B)$.

The multiplication morphism $B\otimes B\rightarrow B$ in ${\rm
Coalg}(B)$ induces a homomorphism in ${\rm Coalg}(A)$
\begin{align}
D:=f^{*}(B\otimes B)\rightarrow C :=f^{*}B.\label{counit}
\end{align}

By  Proposition 2 and the flat cover condition the category $({\rm
Bicomod}_{A}(C), \Box^{C}_{A},C)$ is abelian monoidal ({\em
regularity} condition).

 Since $f$ is a flat cover
the canonical composite
\begin{align}
B ={\rm Hom}_{B}(B\otimes B,B)_{B} \stackrel{f^{*}}{\rightarrow}
{\rm Hom}^{f^{*}B}_{A}(f^{*}(B\otimes B),f^{*}B)^{f^{*}B}_{A} =
{\rm Hom}^{C}_{A}(D,C)^{C}_{A}.\label{inv}
\end{align}
is, by Proposition 2,  an isomomorphism of rings (first part of
the {\em invariants} condition). The isomorphism
\begin{align}f^{*}((B\otimes B)\otimes_{B}(B\otimes B))
\stackrel{\cong}{\rightarrow}f^{*}(B\otimes
B)\square^{f^{*}B}_{A}f^{*}(B\otimes B)
\end{align}
can be rewritten, using $f^{*}(-)=(-)_{B^{o}\otimes B}D$,  as
\begin{align}B\otimes D
\stackrel{\cong}{\rightarrow}D\square^{C}_{A}D\label{inv''}
\end{align}
(second part of the {\em invariants} condition).

By (\ref{inv''}) we have, for $R:=B^{o}\otimes B$ and all $N_{1},
N_{2}\in {\rm Bimod}(B)\simeq\ _{R}{\rm Mod}\simeq {\rm Mod}_{R}$,
the canonical isomorphism
\begin{align}(N_{1}\otimes_{B}N_{2})\otimes_{R}D=
N_{1}\otimes_{R}(B\otimes D
)\otimes_{R}N_{2}\stackrel{\cong}{\rightarrow}
N_{1}\otimes_{R}(D\square^{C}_{A}D)\otimes_{R}N_{2},\label{secan}
\end{align}
so the canonical map (see Lemma \ref{lemma1}) \begin{align}
N_{1}\otimes_{R}(D\square^{C}_{A}D)\otimes_{R} N_{2} \rightarrow
(N_{1}\otimes_{R} D)\square_{A}^{C}( D\otimes_{R} N_{2})
\end{align}
is isomorphic to the canonical homomorphism
\begin{align}f^{*}(N_{1}\otimes_{B}N_{2})\rightarrow
f^{*}N_{1}\square^{f^{*}B}f^{*}N_{2},\label{canhom}
\end{align} which is an isomorphism by flatness of $f$. This
together with Lemma \ref{lemma1} proves the \emph{purity}
condition.

For any $N\in {\rm Bimod}(B)$ we define
\begin{align}
M & :=f^{*}N\in {\rm Bicomod}_{A}(f^{*}B)={\rm
Bicomod}_{A}(C).\nonumber
\end{align}

As in (\ref{inv}) we obtain the following isomorphism of
$B$-bimodules, defined as the canonical composite
\begin{align}
N ={\rm Hom}_{B}(B\otimes B,N)_{B} \stackrel{f^{*}}{\rightarrow}
{\rm Hom}^{f^{*}B}_{A}(f^{*}(B\otimes B),f^{*}N)^{f^{*}B}_{A} =
{\rm Hom}^{C}_{A}(D,M)^{C}_{A},\label{bimod}
\end{align}
where the structure of a $B$-bimodule on the right hand side is
induced by (\ref{inv}) and the canonical structure of ${\rm
Hom}^{C}_{A}(D,C)^{C}_{A}$-bimodule is given as the pair of
canonical composites using the comultiplication $D\rightarrow
D\Box^{C}_{A}D$ and the isomorphisms $M\cong M\Box^{C}_{A}C$,
$M\cong C\Box^{C}_{A}M$
\begin{align}
{\rm Hom}^{C}_{A}(D,M)^{C}_{A}\otimes{\rm
Hom}^{C}_{A}(D,C)^{C}_{A}\rightarrow {\rm
Hom}^{C}_{A}(D\Box^{C}_{A}D,M\Box^{C}_{A}C)^{C}_{A}\rightarrow{\rm
Hom}^{C}_{A}(D,M)^{C}_{A},\nonumber
\end{align}
\begin{align}
{\rm Hom}^{C}_{A}(D,C)^{C}_{A}\otimes{\rm
Hom}^{C}_{A}(D,M)^{C}_{A}\rightarrow {\rm
Hom}^{C}_{A}(D\Box^{C}_{A}D,C\Box^{C}_{A}M)^{C}_{A}\rightarrow{\rm
Hom}^{C}_{A}(D,M)^{C}_{A}.\nonumber
\end{align}
We use (\ref{bimod}) in the following canonical composite
\begin{align}
{\rm Hom}_{B}(B,N)_{B}&
\stackrel{(\ref{bimod})}{\rightarrow} {\rm Hom}_{B}(B,{\rm Hom}^{C}_{A}(D,M)^{C}_{A})_{B}\label{comp}\\
 & ={\rm Hom}^{C}_{A}(B\otimes_{B^{o}\otimes B}D,M)^{C}_{A}\nonumber\\
& = {\rm Hom}^{C}_{A}(D/[D,B],M)^{C}_{A}.\nonumber
\end{align}
On the other hand, we have the following composite similar to a
part of (\ref{bimod})
\begin{align}
{\rm Hom}_{B}(B,N)_{B} \stackrel{f^{*}}{\rightarrow} {\rm
Hom}^{f^{*}B}_{A}(f^{*}B,f^{*}N)^{f^{*}B}_{A} = {\rm
Hom}^{C}_{A}(C,M)^{C}_{A},\label{part}
\end{align}
Since (\ref{comp}) and (\ref{part}) are natural isomorphisms and
the functor $f^{*}: {\rm Bimod}(B)\rightarrow {\rm
Bicomod}_{A}(f^{*}B)$=${\rm Bicomod}_{A}(C)$ is essentially
surjective on objects the canonical morphism $D/[D,B]\rightarrow
C$ in ${\rm Coalg}(A)$ is an isomorphism ({\em freeness}
condition).

Finally, we have the following two natural composites which are
isomorphisms for every $L\in {\rm Bimod}(A)$ and $M$, $N$ as above
\begin{align}
{\rm Hom}_{B}(N,f_{*}L)_{B} & \stackrel{f^{*}}{\rightarrow} {\rm
Hom}^{f^{*}B}_{A}(f^{*}N,f^{*}f_{*}L)^{f^{*}B}_{A} \label{part'}\\
& = {\rm Hom}^{C}_{A}(M,{\rm
Hom}_{A}(D,L)_{A}\otimes_{B^{o}\otimes B }D)^{C}_{A},\nonumber
\end{align}
\begin{align}
{\rm Hom}_{B}(N,f_{*}L)_{B} & = {\rm
Hom}_{A}(f^{*}N,L)_{A}\label{part''}\\
 & = {\rm Hom}_{A}(M,L)_{A}={\rm
Hom}^{C}_{A}(M,C\otimes_{A}L\otimes_{A}C)^{C}_{A}.\nonumber
\end{align}
Similarly, since (\ref{part'}) and (\ref{part''}) are natural
isomorphisms and the functor $f^{*}: {\rm Bimod}(B)\rightarrow
{\rm Bicomod}_{A}(f^{*}B)$=${\rm Bicomod}_{A}(C)$ is essentially
surjective on objects we obtain a natural isomorphism of functors
${\rm Bimod}(A)\rightarrow {\rm Bicomod}_{A}(C)$ ({\em comonad}
condition)
\begin{align}
{\rm Hom}_{A}(D,-)_{A}\otimes_{B^{o}\otimes B }D\rightarrow
C\otimes_{A}-\otimes_{A}C.
\end{align}

({\em Galois extensions $\rightsquigarrow$ Flat covers }) Let us
consider a noncommutative Galois extension ($B\rightarrow {\rm
Hom}_{A}(D,A)_{A}$, $D\rightarrow C$).

By Theorem 1 it determines a morphism $f:{\rm Spec}(A)\rightarrow
{\rm Spec}(B)$, where $f_{*}={\rm Hom}_{A}(D,-)_{A}$,
$f^{*}=-\otimes_{B^{o}\otimes B}D$. By the {\em freeness}
condition and the first part of the {\em invariants} condition we
have the canonical isomorphism in ${\rm Coalg}(A)$
\begin{align}
f^{*}B=B\otimes_{B^{o}\otimes B}D=D/[D,B]\cong C.\label{free}
\end{align}
By (\ref{free}) and the {\em comonad} condition we have an
isomorphism of functors
\begin{align}
f^{*}f_{*}\rightarrow
C\otimes_{A}-\otimes_{A}C=f^{*}B\otimes_{A}-\otimes_{A}f^{*}B.\label{comonad}
\end{align}
In fact, it is an isomorphism of comonads on ${\rm Bimod}(A)$.
Therefore the category of comodules over the comonad $f^{*}f_{*}$
is equivalent to ${\rm Bicomod}_{A}(f^{*}B)$. By the {\em faithful
flatness} condition $f^{*}$ is faithful and exact. Therefore by
the Barr-Beck theorem and the {\em comonad} condition $f^{*}$
induces an equivalence ${\rm Bimod}(B)\rightarrow{\rm
Bicomod}_{A}(f^{*}B)$. By the \emph{regularity} condition the
latter category is monoidal with respect to the cotensor product
over $C=f^{*}B$ with the monoidal unit $C$ and by Proposition 2
this equivalence induced by  $f^{*}$ is (lax) comonoidal. By the
second part of the \emph{invariants} condition and next by the
\emph{purity} condition and Lemma \ref{lemma1} we obtain
isomorphisms in the following canonical decomposition of the
canonical homomorphism (\ref{canhom}) (see also (\ref{secan}))
\begin{align} & f^{*}(N_{1}\otimes_{B}N_{2})  =(N_{1}\otimes_{B}N_{2})\otimes_{R}D \\
 =  N_{1}\otimes_{R}(B\otimes D
)\otimes_{R}N_{2} & \stackrel{\cong}{\rightarrow}
N_{1}\otimes_{R}(D\square^{C}_{A}D)\otimes_{R}N_{2}\stackrel{\cong}{\rightarrow} (N_{1}\otimes_{R} D)\square_{A}^{C}( D\otimes_{R} N_{2})\nonumber\\
 = & f^{*}N_{1}\square^{f^{*}B}f^{*}N_{2}.\nonumber
\end{align}
 This implies that the above equivalence is strong comonoidal, hence $f$ is a flat cover.\
$\blacksquare$

\subsection{Infinitesimals and differential operators} In this
section we show how to define and study infinitesimal structure of
noncommutative schemes. First, we fix the terminology related to
towers and filtrations in abelian categories.

With every tower descending from ${\mathcal F}$
$${\mathcal F}={\mathcal
F}^{0}\twoheadrightarrow\cdots\twoheadrightarrow {\mathcal
F}^{p}\twoheadrightarrow {\mathcal F}^{p+1}\twoheadrightarrow
\cdots \twoheadrightarrow 0$$ we associate an increasing
filtration in ${\mathcal F}$
$$0\rightarrowtail{\mathcal
F}_{0}\rightarrowtail\cdots\rightarrowtail {\mathcal
F}_{p}\rightarrowtail {\mathcal F}_{p+1}\rightarrowtail \cdots
\rightarrowtail {\mathcal F},$$ taking ${\mathcal F}_{p}:={\rm
ker}({\mathcal F}\twoheadrightarrow{\mathcal F}^{p+1})$.

In the dual manner, with every decreasing filtration in ${\mathcal
G}$
$${\mathcal G}={\mathcal G}^{0}\leftarrowtail\cdots\leftarrowtail
{\mathcal G}^{p}\leftarrowtail {\mathcal G}^{p+1}\leftarrowtail
\cdots \leftarrowtail 0$$ we associate a tower ascending to
${\mathcal G}$
$$0\twoheadleftarrow{\mathcal G}_{0}\twoheadleftarrow\cdots\twoheadleftarrow {\mathcal
G}_{p}\twoheadleftarrow {\mathcal G}_{p+1}\twoheadleftarrow \cdots
\twoheadleftarrow {\mathcal G},$$ taking ${\mathcal G}_{p}:={\rm
coker}({\mathcal G}\leftarrowtail{\mathcal G}^{p+1})$.

\vspace{3mm} Let us consider now a morphism of noncommutative
schemes $f:X\rightarrow Y$. For every ${\mathcal F}\in {\rm
Qcoh}(X)$  we define by induction a tower descending from
${\mathcal F}$
\begin{align}f^{*}f_{*}{\mathcal F}^{p}\rightarrow {\mathcal
F}^{p}\rightarrow{\mathcal F}^{p+1}\rightarrow
0,\label{desctow}\end{align} inducing an increasing filtration
${\mathcal F}_{p}$ on ${\mathcal F}$, and for every  ${\mathcal
G}\in {\rm Qcoh}(Y)$ a decreasing filtration in ${\mathcal G}$
\begin{align}
f_{*}f^{*}{\mathcal G}^{p}\leftarrow {\mathcal
G}^{p}\leftarrow{\mathcal G}^{p+1}\leftarrow
0,\label{decfil}\end{align} inducing a tower ${\mathcal G}^{p}$
ascending to ${\mathcal G}$.

Finally, assuming that ${\rm Qcoh}(X)$ is cocomplete and ${\rm
Qcoh}(Y)$ is complete, we can define ${\mathcal F}_{f}:={\rm
colim}_{p}{\mathcal F}_{p}$ and ${\mathcal G}^{f}:={\rm
lim}_{p}{\mathcal G}_{p}$.

In the next proposition we will denote by ${\rm Diff}^{A/K}(M,N)$
differential operators in the sense of Lunts - Rosenberg
\cite{lunros}, which agree with  the definition  of Grothendieck
\cite{gro} if $A$ is commutative.

\begin{proposition} Let $K$ be a commutative ring, $A$ a
$K$-algebra, ${\rm Qcoh}(X)$ consists of $A$-bimodules symmetric
over $K$, ${\rm Qcoh}(Y)$ consists of symmetric $K$-bimodules and
$f^{*}=A\otimes_{K}(-)$. Let ${\mathcal F}:={\rm Hom}_{K}(M,N)$
for some $M,N\in _{A}{\rm Mod}$. Then ${\mathcal F}_{f}={\rm
Diff}^{A/K}(M,N)$.
\end{proposition}

{\em Proof.} By the definition of the increasing filtration
${\mathcal F}_{p}$ the exact sequence (\ref{desctow}) is
isomorphic to
\begin{align}f^{*}f_{*}({\mathcal F}/{\mathcal F}_{p-1})
\rightarrow {\mathcal F}/{\mathcal F}_{p-1}\rightarrow{\mathcal
F}/{\mathcal F}_{p}\rightarrow 0.\end{align} By exactness in the
second term this implies the canonical isomorphism
\begin{align}{\mathcal F}_{p}/{\mathcal F}_{p-1}={\rm im}(f^{*}f_{*}({\mathcal F}/{\mathcal F}_{p-1})
\rightarrow {\mathcal F}/{\mathcal
F}_{p-1}).\label{defind}\end{align} Since $f^{*}=A\otimes_{K}(-)$
has as the right adjoint $f_{*}={\rm Hom}_{A}(A,-)_{A}={\rm
Z}_{A}(-)$, the center of an $A$-bimodule, the functor ${\rm
im}(f^{*}f_{*}(-)\rightarrow (-))$ is nothing but the functor of
the sub-$A$-bimodule generated by the center of a given
$A$-bimodule. Applying this fact to (\ref{defind}) we obtain the
inductive definition of the $p$-th differential part of an
$A$-bimodule ${\mathcal F}$, symmetric over $K$, according to
\cite{lunros}. In the special case of ${\mathcal F}={\rm
Hom}_{K}(M,N)$ one obtains ${\mathcal F}_{p}={\rm
Diff}^{A/K}_{p}(M,N)$, i.e. differential operators of order $\leq
p$. $\blacksquare$

\vspace{3mm}
\paragraph{\textbf{Remark}} In general, as in \cite{lunros}, provided only $f^{*}$ is strongly
comonoidal, the increasing filtration ${\mathcal F}_{p}$ in
${\mathcal F}$ is monoidal, which means that we have natural
transformations
\begin{align}
{\mathcal F}'_{p'}\otimes_{X}{\mathcal F}''_{p''}\rightarrow
({\mathcal F}'\otimes_{X}{\mathcal F}'')_{p'+p''}.
\end{align}
This implies that for every additive category enriched in ${\rm
Qcoh}(X)$ (e.g. enriched in ${\rm Bimod}(A)$ as in \cite{kell}) a
generalized differential part is well defined, generalizing the
category with differential operators as morphisms.
\begin{proposition} Let $A$ be a commutative ring,
$I$ an ideal in $A$, ${\rm Qcoh}(X)$ and ${\rm Qcoh}(Y)$ consist
of symmetric bimodules  over $A/I$ and $A$, respectively, and
$f_{*}$ is the base forgetting from $A/I$ to $A$. Then ${\mathcal
G}^{f}=\hat{{\mathcal G}}_{I}$ ($I$-adic completion).
\end{proposition}

{\em Proof.} Since the forgetting functor $f_{*}$
 has as the left adjoint $f^{*}{\mathcal
G}=A/I\otimes_{A}{\mathcal G}={\mathcal G}/I{\mathcal G}$ the
exact sequence (\ref{decfil}) is isomorphic to
\begin{align}
{\mathcal G}^{p}/I{\mathcal G}^{p}\leftarrow {\mathcal
G}^{p}\leftarrow{\mathcal G}^{p+1}\leftarrow 0,\end{align} which
implies that ${\mathcal G}^{p+1}=I{\mathcal G}^{p}$, hence by
induction ${\mathcal G}^{p}=I^{p}{\mathcal G}$. Therefore by the
definition of the tower ${\mathcal G}_{p}$ ascending to ${\mathcal
G}$ we have ${\mathcal G}_{p}={\mathcal G}/I^{p}{\mathcal G}$.
Finally, passing to the limit we obtain the $I$-adic completion.
$\blacksquare$

\subsection{Cyclic homology of noncommutative schemes} To define cyclic
homologyof noncommutative schemes we need some additional
structure which cannot be derived from the plain abelian monoidal
structure. We derive the axioms of this additional structure from
the canonical structures of the category of affine noncommutative
schemes. First we observe that in the case of $X={\rm Spec}(A)$ we
have a functor
\begin{align}{\rm Tr}_{X}:=(-)\otimes_{A^{o}\otimes A}A\end{align}
\begin{align}{\rm Tr}_{X}: {\rm Qcoh}(X)\rightarrow {\rm Ab}\end{align}
satisfying the following flip symmetry property
\begin{align}{\rm Tr}_{X}({\mathcal F}_{1}\otimes_{X}{\mathcal F}_{2})
\stackrel{\cong}{\rightarrow}{\rm Tr}_{X}({\mathcal
F}_{2}\otimes_{X}{\mathcal F}_{1}).\end{align}

For any geometric morphism $X={\rm Spec}(A)\rightarrow {\rm
Spec}(B)=Y$, which is equivalent to a representation

\begin{align}B\rightarrow{\rm Hom}_{A}(D,A)_{A}\end{align}
we have the following functors

\begin{align}f_{*}={\rm Hom}_{A}(D,-)_{A},\ \  f^{*}=(-)_{B^{o}\otimes B}D\end{align}
By the tensor-hom adjunction they form an adjoint pair
$f^{*}\dashv f_{*}$, i.e. there is an isomorphism of bifunctors
(first hom adjunction axiom)
\begin{align}{\rm Hom}_{X}(f^{*}{\mathcal G}, {\mathcal F})={\rm Hom}_{Y}({\mathcal G}, f_{*}{\mathcal F}).\end{align}
Since $f_{*}$ is a monoidal functor $f^{*}$ is a comonoidal one.
This means that $f_{*}$ transforms algebras (monoids in an abelian
monoidal category) into algebras, while $f^{*}$ transforms
coalgebras (comonoids in an abelian monoidal category) into
coalgebras. In particular $f_{*}{\mathcal O}_{X}$ is always an
algebra while $f^{*}{\mathcal O}_{Y}$ is always a coalgebra.
Structural morphisms ${\mathcal O}_{Y}\rightarrow f_{*}{\mathcal
O}_{X}$ and equivalent to it $f^{*}{\mathcal O}_{Y}\rightarrow
{\mathcal O}_{X}$ are the unit and the counit for $f_{*}{\mathcal
O}_{X}$ and $f^{*}{\mathcal O}_{Y}$.

For the above geometric morphism between affine noncommutative
schemes we have also another pair of functors
\begin{align}f_{!}:=D\otimes_{A^{o}\otimes A}(-),\ \  f^{!}:={\rm Hom}_{B}(D,-)_{B}.\end{align}
By the tensor-hom adjunction they form an adjoint pair
$f_{!}\dashv f^{!}$, i.e. there is an isomorphism of bifunctors
(second hom adjunction axiom)
\begin{align}{\rm Hom}_{Y}(f_{!}{\mathcal F}, {\mathcal G})={\rm Hom}_{X}({\mathcal F}, f^{!}{\mathcal G}),\end{align}

\begin{lemma} The pair of functors
$(f_{!}, f^{*})$ admits two natural transformations (projection
axiom)
\begin{align}f_{!}{\mathcal F}\otimes_{Y} {\mathcal G}\rightarrow
f_{!}({\mathcal F}\otimes_{X} f^{*}{\mathcal G})\end{align}
\begin{align}{\mathcal G}\otimes_{Y} f_{!}{\mathcal
F}\rightarrow f_{!}( f^{*}{\mathcal G}\otimes_{X}{\mathcal F})
\end{align}
\end{lemma}

\emph{Proof.} Using the equivalence of bimodules and left or right
modules over the enveloping ring it is enough to construct only
the last transformation, the previous one we obtain in an
analogical way, changing the role of left and right modules.

\begin{eqnarray}{\mathcal G}\otimes_{Y} f_{!}{\mathcal
F} & = & {\mathcal G}\otimes_{B} (D\otimes_{A^{o}\otimes A
}{\mathcal F})\\
& = & ((B\otimes B)\otimes_{B}{\mathcal G})\otimes_{B^{o}\otimes
B}D\otimes_{A^{o}\otimes A}{\mathcal
F}\\
& = & (f^{*}((B\otimes B)\otimes_{B}{\mathcal
G})\otimes_{A}{\mathcal F})\otimes_{A^{o}\otimes A}A\\
& \rightarrow & (f^{*}(B\otimes B)\otimes_{A}f^{*}{\mathcal
G}\otimes_{A}{\mathcal F})\otimes_{A^{o}\otimes A}A\\
& = & (D\otimes_{A}f^{*}{\mathcal G}\otimes_{A}{\mathcal
F})\otimes_{A^{o}\otimes A}A\\
& = & D\otimes_{A^{o}\otimes A}(f^{*}{\mathcal
G}\otimes_{A}{\mathcal F})\\
& = & f_{!}( f^{*}{\mathcal G}\otimes_{X}{\mathcal F}).
\end{eqnarray}

$\blacksquare$

\begin{lemma} There exists a natural isomorphism of bifunctors (trace adjunction axiom)
\begin{align}{\rm Tr}_{X}({\mathcal
F}\otimes_{X}f^{*}{\mathcal G})\cong {\rm Tr}_{Y}(f_{!}{\mathcal
F}\otimes_{Y}{\mathcal G})\end{align}

\end{lemma}

\emph{Proof.} Using the flip transformations under ${\rm Tr}_{X}$
and ${\rm Tr}_{Y}$ it is enough to construct the flip-equivalent
transformation.

\begin{eqnarray}{\rm Tr}_{X}(f^{*}{\mathcal G}\otimes_{X}{\mathcal
F}) & = & (({\mathcal G}\otimes_{B^{o}\otimes
B}D)\otimes_{A}{\mathcal
F})\otimes_{A^{o}\otimes A}A\\
& = & {\mathcal G}\otimes_{B^{o}\otimes B}D\otimes_{A^{o}\otimes
A}{\mathcal F}\\
& = & ({\mathcal G}\otimes_{B}(D\otimes_{A^{o}\otimes A}{\mathcal
F}))\otimes_{B^{o}\otimes B}B\\
& = & {\rm Tr}_{Y}({\mathcal G}\otimes_{Y}f_{!}{\mathcal
F}).\end{eqnarray} $\blacksquare$

\vspace{3mm}
\paragraph{\textbf{Definition.}} Let $p:X\rightarrow S$ be a morphism of noncommutative
schemes. A \emph{system of coefficients (relative to $S$)} is an
object $M\in {\rm Qcoh}(X)$ equipped with a braiding
$$\beta: p^{*}{\mathcal O}_{S}\otimes_{X} M\rightarrow M\otimes_{X} p^{*}{\mathcal O}_{S},$$
such that the following diagrams commute (where $C:=p^{*}{\mathcal
O}_{S}\in {\rm Coalg}(X)$, $\otimes :=\otimes_{X}$, transpositions
of $M$ and copies of $C$ are obtained via $\beta$ and
$\beta^{-1}$, and $\Delta_{i}$, $\varepsilon_{i}$ come from
$\Delta : C\rightarrow C \otimes C$ and $\varepsilon :
C\rightarrow {\mathcal O}_{X}$ applied to the $i$-th factor.)

\vspace{3mm}\paragraph{\emph{Diagram I.}}
$$\begin{array}{ccc}
C\otimes C\otimes M & \rightarrow & M\otimes C\otimes C\\
\Delta_{2}\downarrow & & \downarrow \Delta_{2}\\
C\otimes C\otimes C\otimes M & \rightarrow & M\otimes C\otimes
C\otimes C\end{array}
$$
\emph{Diagram II.}
$$\begin{array}{ccccccc}
M\otimes C  & \stackrel{\Delta_{2}}{\rightarrow} & M\otimes
C\otimes C & \rightarrow & C\otimes M\otimes C &
\stackrel{\Delta_{1}}{\rightarrow} & C\otimes C\otimes M\otimes C \\
\downarrow & & & & & & \downarrow \\
C\otimes M & \stackrel{\Delta_{1}}{\rightarrow} & C\otimes
C\otimes M & \rightarrow & C\otimes M\otimes C&
\!\!\stackrel{\Delta_{3}}{\rightarrow}\!\! & C\otimes M\otimes
C\otimes C
\end{array}
$$
\emph{Diagram III.}
$$\begin{array}{ccccccc}
C\otimes M\otimes C  & \stackrel{\Delta_{1}}{\rightarrow} &
C\otimes C\otimes M\otimes C & \rightarrow & C\otimes M\otimes
C\otimes C & \stackrel{\varepsilon_{3}}{\rightarrow} & C\otimes
M\otimes C \\
\downarrow & & & & & & \downarrow \\
M\otimes C\otimes C  & \stackrel{\varepsilon_{2}}{\rightarrow} &
M\otimes C & \rightarrow & C\otimes M&
\stackrel{\Delta_{1}}{\rightarrow} & C\otimes C\otimes M
\end{array}
$$
\emph{Diagram IV.}
$$\begin{array}{ccc}
C\otimes M & \rightarrow & M\otimes C\\
\Delta_{1}\downarrow & & \uparrow \varepsilon_{2}\\
C\otimes C\otimes M & \rightarrow & M\otimes C\otimes C\end{array}
$$
\emph{Diagram V.}
$$\begin{array}{ccc}
C\otimes C\otimes M  & \rightarrow & C\otimes M\otimes C \\
\varepsilon_{2}\downarrow & & \downarrow\varepsilon_{3}\\
C\otimes M& & C\otimes M\\
\downarrow & & \downarrow
\\
M\otimes C & & M\otimes C \\
\ \ \ \ \ \ \ \ \ \ \ \ \ \varepsilon_{2}\searrow & & \swarrow\varepsilon_{2}\ \ \ \ \ \ \ \ \ \ \ \ \   \\
 & M &
\end{array}
$$
\emph{Diagram VI.}
$$\begin{array}{ccc}
C\otimes C\otimes M  & \rightarrow & M\otimes C\otimes C \\
\varepsilon_{2}\downarrow & & \downarrow\varepsilon_{2}\\
C\otimes M&\rightarrow & M\otimes C
\end{array}
$$

\vspace{3mm}
\paragraph{\textbf{Definition.}} $M$ is a
\textsl{cyclic system of coefficients} relative to $S$ if for all
$p,q>0$ the following diagrams commute (braiding-trace
compatibility):

$$\begin{array}{ccc}
{\rm Tr}_{X}(C^{\otimes p}\otimes M\otimes C^{\otimes q}) & \rightarrow & {\rm Tr}_{X}(C^{\otimes p+1}\otimes M\otimes C^{\otimes q-1})\\
\downarrow & & \downarrow \\
{\rm Tr}_{X}(C^{\otimes p-1}\otimes M\otimes C^{\otimes q+1}) &
\rightarrow & {\rm Tr}_{X}(C^{\otimes p}\otimes M\otimes
C^{\otimes q}),\end{array}$$

$$\begin{array}{ccc}
{\rm Tr}_{X}(C^{\otimes p}\otimes M) & \rightarrow & {\rm Tr}_{X}(M\otimes C^{\otimes p})\\
\downarrow & & \downarrow \\
{\rm Tr}_{X}(C^{\otimes p-1}\otimes M\otimes C) & \rightarrow &
{\rm Tr}_{X}(C^{\otimes p}\otimes M),\end{array}$$ where vertical
arrows are induced by braiding transpositions $C\otimes M
\rightarrow M\otimes C$ or their inverses,  and horizontal ones
are induced by natural flip isomorphisms of ${\rm Tr}_{X}$.

\vspace{3mm}
\paragraph{\textbf{Example.}} $C=p^{*}{\mathcal O}_{S}$ itself is
a (trivial) cyclic system of coefficients relative to $S$ with the
identity braiding $\beta: C\otimes C\rightarrow C\otimes C$.

\vspace{3mm}
\paragraph{\textbf{Definition.}} A \emph{cyclic object}
in an abelian category consists of a collection of morphisms
$\partial_{i}:C_{n}\rightarrow C_{n-1}$, $s_{i}:C_{n}\rightarrow
C_{n+1}$, $i=0,\ldots,n$, and $t_{n}:C_{n}\rightarrow C_{n}$,
satisfying
\begin{eqnarray}\partial_{i}\partial_{j} & = & \partial_{j-1}\partial_{i},  \ \ \ \ \ \ \ \  i<j,\\
                              s_{i}s_{j} & = & s_{j+1}s_{i},  \ \ \ \ \ \ \ \  i\leq j,\\
                       \partial_{i}s_{j} & = &  \left\{ \begin{array}{lll}
                               s_{j-1}\partial_{i}, & \ \   i<j,\\
                               {\rm id}, & \ \   i=j, j+1,\\
                               s_{j}\partial_{i-1}, & \ \    i>j+1,
                               \end{array}\right.\\
\partial_{i}t_{n} & = & t_{n-1}\partial_{i-1}, \ \ \ \ \ i=1,\ldots ,n,\\
\partial_{0}t_{n} & = & \partial_{n},  \\
s_{i}t_{n} & = & t_{n+1}s_{i-1},  \ \ \ \ \ i=1,\ldots ,n,\\
s_{0}t_{n} & = & t_{n+1}^{2}s_{n},  \\
t_{n}^{n+1} & = & 1.  \end{eqnarray}

\vspace{3mm}
\paragraph{\textbf{Definition.}} For every cyclic system $M$ of coefficients relative to $S$ we define
$$C_{X/S}(M)_{n}:={\rm Tr}_{X}(M\otimes_{X} (p^{*}{\mathcal
O}_{S})^{\otimes_{X} n}),$$ and  (in the element-wise convention!)

$$\partial_{i}:C_{X/S}(M)_{n}\rightarrow C_{X/S}(M)_{n-1},$$
\begin{align*}
\partial_{0}{\rm Tr}_{X}(m,c_{1},\ldots,c_{n}) & :=  {\rm Tr}_{X}(m\varepsilon(c_{1}),c_{2},\ldots,c_{n})  \\
\partial_{i}{\rm Tr}_{X}(m,c_{1},\ldots,c_{n}) & :=  {\rm Tr}_{X}(m,c_{1},\ldots,c_{i}\varepsilon(c_{i+1}),\ldots,c_{n}),   \ i=1,\ldots,
n-1,\\
\partial_{n}{\rm Tr}_{X}(m,c_{1},\ldots,c_{n}) & :=  {\rm Tr}_{X}((1\otimes \varepsilon )\beta (c_{n},m),c_{1},\ldots,c_{n-1})
\end{align*}

$$s_{i}:C_{X/S}(M)_{n}\rightarrow C_{X/S}(M)_{n+1},$$
\begin{align*}
s_{i}{\rm Tr}_{X}(m,c_{1},\ldots,c_{n}) & := {\rm
Tr}_{X}(m,c_{1},\ldots,\Delta(c_{i+1}),\ldots,c_{n}),   \
i=0,\ldots,
n-1,\\
s_{n}{\rm Tr}_{X}(m,c_{1},\ldots,c_{n}) & := \sum{\rm
Tr}_{X}(\beta(c'_{1(2)},m'),c_{2},\ldots,c_{n},c'_{1(1)}),
\end{align*}
where
$$\beta^{-1}(m,c)=\sum c'\otimes m',\ \Delta(c)=\sum c_{(1)}\otimes c_{(2)},$$

$$t_{n}:C_{X/S}(M)_{n}\rightarrow
C_{X/S}(M)_{n},$$
$$t_{n}{\rm Tr}_{X}(m,c_{1},\ldots,c_{n}):={\rm Tr}_{X}(\beta(c_{n}, m), c_{1},\ldots,c_{n-1}).$$

\vspace{3mm}
\paragraph{\textbf{Definition.}} We say that ${\rm Tr}_{X}$ is
\emph{faithful relative to $S$} if for every $n>0$ the functor
$$C_{X/S}(-)_{n}: {\rm Qcoh}(X)\rightarrow {\rm Ab}$$
is faithful.

\begin{theorem} For every cyclic system $M$ of coefficients
relative to $S$ the system $(C_{X/S}(M)_{\bullet},
\partial_{\bullet}, s_{\bullet}, t_{\bullet})$ is a cyclic
object in the category of abelian groups.

Assume that ${\rm Tr}_{X}$ is faithful relative to $S$. Then for
every $M$ with a braiding with respect to $p^{*}{\mathcal O}_{S}$
compatible with ${\rm Tr}_{X}$ the para-cyclic relations between
$(\partial_{\bullet}, s_{\bullet}, t_{\bullet})$ defined as above
define on $M$ a structure of a cyclic system of coefficients
relative to $S$.
\end{theorem}

\emph{Proof}. Diagram V implies  $\partial_{n-1}\partial_{n}  =
\partial_{n-1}\partial_{n-1}$.
 Diagram II
implies $s_{0}s_{n}=s_{n+1}s_{0}$. Coassociativity of $\Delta:
C\rightarrow C\otimes C $ implies $s_{i}s_{i}=s_{i+1}s_{i}$ for
$i=0,\ldots, n-1$. Coassociativity and invertibility of braiding
imply $s_{n}s_{n}=s_{n+1}s_{n}$. Diagram III implies
$\partial_{0}s_{n}  =  s_{n-1}\partial_{0}$. The left counit
property $(\varepsilon\otimes {\rm id})\Delta={\rm id}$ implies
$\partial_{i}s_{i}  =  {\rm id}$ and $\partial_{i+1}s_{i}  = {\rm
id}$ for $i=0,\ldots, n-1$. Diagram IV implies
$\partial_{n+1}s_{n}  = {\rm id}$.  Diagram VI implies
$\partial_{n}t_{n} = t_{n-1}\partial_{n-1}$. Diagram I implies
$s_{0}t_{n}=t_{n+1}^{2}s_{n}$. The trace flip-braiding
compatibility implies $t_{n}^{n+1}={\rm id}$. All other cyclic
object relations are fulfilled automatically by the definition of
$(\partial_{\bullet}, s_{\bullet}, t_{\bullet})$.

If ${\rm Tr}_{X}$ is faithful relative to $S$, all implications
between commutativity of diagrams I-VI and cyclic object relations
become equivalences.

$\blacksquare$

We denote the respective Hochschild, cyclic and periodic cyclic
homology by $HH_{X/S}(M)_{\bullet}$, $HC_{X/S}(M)_{\bullet}$ and
$HP_{X/S}(M)_{\bullet}$, respectively.

\vspace{3mm}
\paragraph{\textbf{Example.}} Let $K\rightarrow A$ be a central
ring homomorphism from a commutative ring $K$, regarded as a
geometric morphism $p:{\rm Spec}(A)=X\rightarrow S={\rm Spec}(K)$
(both spectra noncommutative). Then $C=p^{*}{\mathcal O
}_{S}=A\otimes_{K}A$. We have also ${\rm Tr}_{X}{\mathcal
F}:=A\otimes_{A^{o}\otimes A}{\mathcal F}$ and we can take $M:=C$.
Then we obtain the classical cyclic object of a $K$-algebra $A$.

\subsubsection{Functoriality of cyclic objects.}
Assume now that we have a commutative diagram

$$\begin{array}{ccccc}
         & f         &     \\
X     & \overrightarrow{\ \ \ \ \ \ \ } &  Y, \\
 \ \ \ \ \ \ \ \ \ \ \ \ \        p  \searrow  &        &    \swarrow q  \ \ \ \ \ \ \ \ \ \ \ \ \  \\
         & S &
\end{array}$$
in the category of noncommutative schemes over a noncommutative
scheme $S$ equipped with trace functors and with the geometric
morphism $f^{*}\dashv f_{*}$ completed by a compatible pair
$f_{!}\dashv f^{!}$. Assume that $M$ and $N$ are cyclic systems of
coefficients relative to $S$ on $X$ and $Y$, respectively. Assume
that we have a morphism $N\rightarrow f_{!}M$ making the canonical
diagram
$$\begin{array}{ccc}
q^{*}{\mathcal O}_{S}\otimes_{Y}N & \rightarrow & N
\otimes_{Y}q^{*}{\mathcal O}_{S}\\
\downarrow & & \downarrow\\
f_{!}(p^{*}{\mathcal O}_{S}\otimes_{X}M) & \rightarrow & f_{!}(M
\otimes_{X}p^{*}{\mathcal O}_{S}) \end{array}$$ commutative.

Then it induces a morphism of cyclic objects

\begin{align}(C_{X/S}(M)_{\bullet},
\partial_{\bullet}, s_{\bullet}, t_{\bullet})\leftarrow(C_{Y/S}(N)_{\bullet},
\partial_{\bullet}, s_{\bullet}, t_{\bullet}).\end{align}
and hence the morphisms of homologies
\begin{align}HH_{X/S}(M)_{\bullet}\leftarrow(HH_{Y/S}(N)_{\bullet},\\
HC_{X/S}(M)_{\bullet}\leftarrow(HC_{Y/S}(N)_{\bullet},\\
HP_{X/S}(M)_{\bullet}\leftarrow(HP_{Y/S}(N)_{\bullet}.\end{align}

\subsubsection{Comparison with other constructions.}
\vspace{3mm} In a recent paper \cite{bost} autors consider a
construction of cyclic objects based on comonads and
distributivity laws. Down to the earth, restronged to the context
of our construction, their structure is based on the following two
diagrams

{\emph{\flushleft{Diagram I'}}}
$$\begin{array}{ccc}
C\otimes M & \rightarrow & M\otimes C\\
\Delta_{1}\downarrow & & \downarrow \Delta_{2}\\
C\otimes C\otimes M & \rightarrow & M\otimes C\otimes C\end{array}
$$
\emph{\flushleft{Diagram V'}}
$$\begin{array}{ccc}
C\otimes M &  \rightarrow & M\otimes C \\
\ \ \ \ \ \ \ \ \ \ \ \ \ \varepsilon_{1}\searrow & & \swarrow\varepsilon_{2}\ \ \ \ \ \ \ \ \ \ \ \ \   \\
 & M. &
\end{array}
$$

The first difference between this cyclic object and our consists
in the number of copies of a coalgebra $C$ in every degree in the
complex. In degree $n$ we see in their complex $n+1$ copies while
in our we have only $n$ of them. This suggests that these two
constructions have different flavor. Indeed, the classical
(commutative) object corresponding to our construction is the
DeRham cohomology with values in a module with an integrable
connection. In the case of the Sweedler coalgebra
$C=A\otimes_{K}A$ in ${\rm Bimod}(A)$ for a commutative
$K$-algebra $A$ over a commutative ring $K$ we think of the
bimodule $M$ as of a noncommutative analog of a sheaf supported on
the first infinitesimal neighborhood of the diagonal in the
cartesian square of a scheme over $K$. It contains the information
about its restriction to the diagonal together with an
infinitesimal variation of this restriction encoded in the brading
$\beta$, i. e. a connection in the Grothendieck approach. Our
diagrams I-VI together with  the trace-braiding compatibility form
a noncommutative analog of the property of being supported on the
first infinitesimal neighborhood and integrability of the
connection.

The conditions on coefficients $M$ in the construction of
\cite{bost} is a generalization of the stable anti-Yetter-Drinfeld
condition from Hopf-cyclic homology with coefficients (in the dual
approach of Jara-\c{S}tefan). This condition means that $M$ is
regarded as a noncommutative analog of a stable equivariant sheaf
on the group of symmetries acted by itself via conjugations. This
means that it depends only on symmetries of the space and not on
the space itself. Therefore it requires only diagrams I' and VI',
which play the role of our diagrams I and VI. Moreover, the
passage to the cyclic object consists there in quotiening the
paracyclic object by the relation forcing the cyclic relation
$t^{n+1}_{n}={\rm id}$, a la Kaygun. It is an analog of dividing
the DeRham complex tensored by a module with a non-integrable
connection by the image of the curvature to obtain a complex, and
hence does not correspond to integrability of the connection.

In another recent paper \cite{kal} the author uses a formalism of
cocartesian objects in symmetric monoidal categories. In the
simplest case, an algebra $A$ over a commutative ring $K$ is
considered, and we can pass to our context taking the respective
Sweedler construction, i.e.
$$C=A\otimes_{K}A.$$
The author assumes that there is given a twist
$$A\otimes_{K}M\rightarrow M\otimes_{K}A$$
in the category of bimodules over the algebra $A\otimes_{K}A$,
satisfying some cocycle condition. It can be compared with our
braiding, taking into account the fact that it can be written as
$$C\otimes_{A}M=A\otimes_{K}M\rightarrow M\otimes_{K}A=M\otimes_{A}C.$$
The respective commutative diagram encoding this cocycle condition
is of the form

$$\begin{array}{ccc}
A\otimes_{K} A\otimes_{K} M &  \longrightarrow & A\otimes_{K} M\otimes_{K} A \\
\ \ \ \ \ \ \ \ \ \ \ \ \ \searrow & & \swarrow\ \ \ \ \ \ \ \ \ \ \ \ \   \\
 & M\otimes_{K} A\otimes_{K} A &
\end{array}
$$
where the south-east arrow is defined as the transposition of the
first and the third factor. It is clear that one has to use the
symmetry of the monoidal category to do this. Also the algebra
structure on tensor powers of $A$ over $K$, e.g. on $A\otimes_{K}
A$, needs this symmetry. Therefore this construction makes no
sense over a noncommutative base ring $K$. Moreover, expressing
the latter diagram in terms of the Sweedler construction one gets
$$\begin{array}{ccc}
C\otimes_{A} C\otimes_{A} M &  \longrightarrow & C\otimes_{A} M\otimes_{A} C \\
\ \ \ \ \ \ \ \ \ \ \ \ \ \searrow & & \swarrow\ \ \ \ \ \ \ \ \ \ \ \ \   \\
 & M\otimes_{A} C\otimes_{A} C &
\end{array}
$$
but the south-east arrow cannot be defined purily in terms of the
category of $A$-bimodules, without referring to the special
structure of $C$.


\begin{thebibliography}{99}

\bibitem{art} Artin, M.; Zhang, J. J.: \emph{Noncommutative projective schemes}, Adv.
Math. 109 (1994), no. 2, pp. 228-287.

\bibitem{bost} B\"{o}hm, G.; \c{S}tefan, D.: \emph{(Co)cyclic (co)homology of
bialgebroids: An approach via (co)monads.} arXiv:0705.3190v1
[math.KT] 22 May 2007


\bibitem{brz} Brzezi\'{n}ski, T.: {\em On modules associated to coalgebra-Galois
extensions.} J. Algebra 215 (1999), 290--317.



\bibitem{brz'} Brzezi\'{n}ski, T.: {\em The structure of corings. Induction functors,
Maschke-type theorem, and Frobenius and Galois-type properties.}
Alg. Rep. Theory 5 (2002), 389--410.

\bibitem{brz''} Brzezi\'{n}ski, T.: {\em Galois comodules.}
J. Algebra 290 (2005), 503--537.

\bibitem{brzhaj} Brzezi\'{n}ski, T.; Hajac, P.M.: {\em Coalgebra extensions and algebra coextensions of Galois type.}
Comm. Algebra 27 (1999), 1347--1367.

\bibitem{brzwis} Brzezi\'{n}ski, T.; Wisbauer, R.: {\em Galois comodules.}
London Math. Soc. LNS. 309, Cambridge University Press, 2003.

\bibitem{bal} Balmer, P.: {\em The spectrum of prime ideals in
tensor differential graded categories.} J. Reine Angew. Math. 588
(2005), 149--168.





\bibitem{ber}  Bergman, G. M.; Hausknecht, A. O.: {\em Cogroups
and co-rings in categories of associative rings,} Mathematical
Surveys and Monographs, vol. 45, American Mathematical Society,
Providence, RI, 1996

\bibitem{cae} Caenepeel, S.: {\em Galois corings from the descent theory point of
view.} Galois theory, Hopf algebras, and semiabelian categories,
163--186, Fields Inst. Commun., 43, Amer. Math. Soc., Providence,
RI, 2004.

\bibitem{caegrover} Caenepeel, S.; De Groot, E.; Vercruysse, J.:
{\em Galois theory for comatrix corings: Descent theory, Morita
theory, Frobenius and separability properties.} Trans. Amer. Math.
Soc., ({\em to appear}), arXiv:math.RA/0406436 (2004).

\bibitem{car} Cartier, P.: {\em A mad day's work:
from Grothendieck to Connes and Kontsevich - the evolution of
concepts of space and symmetry.} Bull. AMS 38 (2001), no. 4,
389--408.


\bibitem{eil} Eilenberg, S.: \emph{Abstract description of some basic functors}, J.
Indian Math. Soc. (N.S.) 24 (1960), 231-234.




\bibitem{elkgom} El Kaoutit, L.; G\'omez-Torecillas, J.: {\em Comatrix corings:
Galois corings, descent theory, and a structure theorem for
cosemisimple corings.} Math. Z.  244 (2003), 887--906.

\bibitem{fre} Freyd, P.: {\em Algebra valued functors in general and tensor products in
particular.} Colloquium Mathematicum (Wroc\l aw), 14 (1966),
89--106.

\bibitem{gab} Gabriel, P.: Des cat\'egories ab\'eliennes, Bull. Soc. Math.
France 90 (1962), 323-448.

\bibitem{gol} Golan,J. S.; Raynaud,J.; van Oystaeyen, F.: Sheaves over the
spectra of certain noncom- mutative rings, Comm. Alg. 4(5),
(1976), pp. 491-502.

\bibitem{gom} G\'omez-Torrecillas, J.: \emph{Separable functors in corings.} Int. J. Math.
Math. Sci. 30 (2002), no. 4, 203--225.

\bibitem{gro} Grothendieck, A.: {\em \'El\'ements de g\'eom\'etrie alg\'ebrique. \'Etude locale de sch\'emas.}
Publ. Math. IHES 32 (1967).

\bibitem{gro'} Grothendieck, A. et al.: {\em R\^evetements \'etale et group fondamental.}
S\'eminaire de G\'eometrie Alg\'ebrique du Bois Marie 1960-1961
(SGA 1), LNM 224, Springer, 1971.

\bibitem{gru} Grunenfelder, L.;  Par\'e, R.: \emph{Families parametrized by coalgebras},
J. Algebra, 107 (1987), 316-375.


\bibitem{hod} Hodges, T. J.; Smith, S. P.: \emph{Sheaves of noncommutative algebras
and the Beilinson- Bernstein equivalence of categories}, Proc.
AMS, 92, N. 3, 1985.

\bibitem{kal} Kaledin, D.: {\em  Cyclic homology with coefficients}, arXiv:math.KT/0702068v1 3 Feb 2007


\bibitem{kell} Kelly, G. M.: {\em  Basic Concepts of Enriched Category Theory}, London
Math. Soc. Lecture Notes Series 64 (Cambridge University Press
1982).

\bibitem{kon} Kontsevich, M.: {\em Triangulated categories and geometry}, Course at
the \'{E}cole Normale Sup\'{e}rieure, Paris, Notes taken by J. Bella\"{\i}che,
J.-F. Dat, I. Marin, G. Racinet and H. Randriambololona, 1998.

\bibitem{konros} Kontsevich, M., Rosenberg, A. L.: {\em Noncommutative smooth spaces.}
The Gelfand Math. Seminars, 1996--1999, pp. 85--108, Birkh\"{a}user,
2000.

\bibitem{konros'} Kontsevich, M., Rosenberg, A. L.: {\em Noncommutative
spaces,} MPI 2004 - 35; {\em Noncommutative spaces and flat
descent,} MPI 2004 - 36; {\em Noncommutative stacks,} MPI 2004 -
37; preprints, Bonn 2004.

\bibitem{lunros} Lunts, V., Rosenberg, A. L.: {\em Differential operators on noncommutative
rings,} Sel. Math., New Ser. 3 (1997), 335-359.

\bibitem{lunros'} Lunts, V., Rosenberg, A. L.: {\em Differential calculus in noncommutative
algebraic geometry I,} Preprint MPI 96-53.

\bibitem{lunros''} Lunts, V., Rosenberg, A. L.: {\em Differential calculus in noncommutative
algebraic geometry II,} Preprint MPI 96-76.

\bibitem{lunros'''} Lunts, V. A., Rosenberg, A. L.: \emph{Localization for quantum
groups}, Selecta Math. (N.S.) 5 (1999), no. 1, pp. 123-159.

\bibitem{mac} Mac Lane, S.: {\em Natural associativity and commutativity}, Rice Univ.
Stud. 49 (1963), 28--46.


\bibitem{macmoe} Mac Lane, S., Moerdijk, I.: \emph{Sheaves in geometry and logic},
Springer 1992.

\bibitem{mur} Murdoch, D. C.; van Oystaeyen, F.: \emph{Noncommutative localization and
sheaves}, J. Alg. 38 (1975), pp. 500-515.

\bibitem{nus} Nuss, P.: \emph{Noncommutative descent and non-abelian cohomology}.
Algebra, 3. J. Math. Sci. 82 (1996), no. 6, 3824-3831.

\bibitem{orl} Orlov, D. O.: \emph{Quasi-coherent sheaves in commutative and
non-commutative geometry}, Izvestiya: Math. 67:3 (2003), pp.
119-138.

\bibitem{hai} Phung Ho Hai: {\em An embedding theorem for abelian monoidal categories}, Compositio
Math., 132(2) (2002), 27--48.

\bibitem{ros} Rosenberg, A. L.: \emph{The spectrum of abelian categories and
reconstruction of schemes}, in \emph{Rings, Hopf algebras, and
Brauer groups (Antwerp/Brussels, 1996)}, pp. 257-274, Lec. Notes
Pure Appl. Math. 197, Dekker, NY 1998.

\bibitem{ros'} Rosenberg, A. L.: \emph{Underlying spectra of noncommutative schemes},
MPI-2003-111.

\bibitem{schsch} Schauenburg, P.; Schneider, H.-J.: {\em On generalized Hopf Galois extensions} arXiv:math.QA/0405184 (2004).

\bibitem{schn} Schneider, H.-J.: {\em Principal homogeneous spaces for arbitrary Hopf algebras}, Israel J. Math., 72 (1-2) (1990), 167-195.  arXiv:math.QA/0405184 (2004).

\bibitem{tak''} Takeuchi, M.: \emph{A note on geometrically reductive groups}, J. Fac. Sci., Univ. Tokyo, Sect. 1, 20 (3) (1973),
384-396.


\bibitem{tak} Takeuchi, M.: \emph{Introduction to $\sqrt{Morita}$
theory}. Proceedings of the 17th symposium on ring theory
(Tsukuba, 1984), 78-86, Okayama Univ., Okayama, 1984.

\bibitem{tak'} Takeuchi, M.: \emph{$\sqrt{Morita}$
theory-formal ring laws and monoidal equivalences of categories of
bimodules}, J. Math. Soc. Japan 39 (1987), no. 2, 301-336.

\bibitem{ver} Verschoren, A.: \emph{Sheaves and localization}, J. Algebra 182 (1996),
no. 2, pp. 341-346.

\bibitem{wat} Watts, C. E.: \emph{Intrinsic characterization of some additive
functors}, Proc. AMS 11 (1960), pp. 1-8.

\bibitem{wis} Wisbauer, R.: {\em On Galois corings.} Hopf algebras in non-commutative geometry and physics,
S. Caenepeel and F. van Oystaeyen (eds.), LNPAM 239, Marcel
Dekker, New York, 2004.

\bibitem{wis'} Wisbauer, R.: {\em On Galois comodules.} Comm. Algebra 34 (2006), 2683-2711.

\end{thebibliography}
\end{document}